\documentclass[12pt]{article} 
\usepackage[sectionbib]{natbib}
\usepackage[caption=false]{subfig}
\usepackage{array,epsfig,fancyheadings,rotating}
\usepackage[]{hyperref}  

\usepackage{amsmath}
\usepackage{amssymb}
\usepackage{amsfonts}
\usepackage{multirow}
\usepackage{amsthm}

\newtheorem{theorem}{Theorem}

\newtheorem{corollary}{Corollary}

\theoremstyle{definition}

\newtheorem{remark}{Remark}
\begin{document}
	
	
	\renewcommand{\baselinestretch}{2}
	
	\markright{ \hbox{\footnotesize\rm Statistica Sinica
		}\hfill\\[-13pt]
		\hbox{\footnotesize\rm
		}\hfill }
	
	\markboth{\hfill{\footnotesize\rm Abhik Ghosh AND Ayanendranath Basu} \hfill}
	{\hfill {\footnotesize\rm Robust Bounded Influence Tests for I-NH Obs.} \hfill}
	
	\renewcommand{\thefootnote}{}
	$\ $\par
	
	
	\fontsize{12}{14pt plus.8pt minus .6pt}\selectfont \vspace{0.8pc}
	\centerline{\large\bf Robust Bounded Influence Tests for }
	\vspace{2pt} \centerline{\large\bf Independent Non-Homogeneous Observations
}
	\vspace{.4cm} \centerline{Abhik Ghosh and Ayanendranath Basu} \vspace{.4cm} \centerline{\it
		Indian Statistical Institute, Kolkata, India} \vspace{.55cm} \fontsize{9}{11.5pt plus.8pt minus
		.6pt}\selectfont
	

\begin{quotation}
	\noindent {\it Abstract:}
Experiments often yield non-identically distributed data 
for statistical analysis.
Tests of hypothesis under such set-ups are generally performed using 
the likelihood ratio test, which is non-robust with respect to outliers and model misspecification. 
In this paper, we consider the set-up of non-identically but independently distributed observations
and develop a general class of test statistics for testing parametric hypothesis 
based on the density power divergence. 
The proposed tests have bounded influence functions, are highly robust with respect to data contamination,
have high power against contiguous alternatives, and are consistent at any fixed alternative.
The methodology is illustrated by the simple and generalized linear regression models with fixed covariates.
\par
\vspace{9pt}
	\par
\end{quotation}\par

\def\thefigure{\arabic{figure}}
\def\thetable{\arabic{table}}

\renewcommand{\theequation}{\thesection.\arabic{equation}}

\fontsize{12}{14pt plus.8pt minus .6pt}\selectfont

\section{Introduction}\label{SEC:intro}

One of the most important paradigms of statistical inference is hypothesis testing;
arguably the most common test is the likelihood ratio test (LRT).
However, like the maximum likelihood estimator (MLE), 
the LRT may lead to unstable inference in the presence of outliers.
Attempts to rectify this (\cite{Simpson:1989, Lindsay:1994, Basu/etc:2013a, Basu/etc:2013b})
have mostly been in the context of independent and identically distributed (i.i.d.) data. 
In this paper, we consider the general case of non-identically distributed data. 
Suppose the observed data $Y_1, \ldots, Y_n$ are independent, but 
for each $i$, $Y_i \sim g_i$ with $g_1, \ldots, g_n$ being possibly different densities with respect 
to some common dominating measure. We model $g_i$ by the family 
${\mathcal F}_{i, {\boldsymbol{\theta}}} = \{f_i(\cdot;{\boldsymbol{\theta}}) |~ {\boldsymbol{\theta}} \in \Theta \}$ for all $i=1,2, \ldots, n$. 
Let $G_i$ and $F_i(\cdot,{\boldsymbol{\theta}})$ be the indicated distribution functions. 
Even though the $Y_i$s  have possibly different densities, 
they share the common  parameter ${\boldsymbol{\theta}}$. 
We will refer to this set-up as the independent non-homogeneous (I-NH) set-up.

The most prominent application of the I-NH set-up is the regression model 
with non-stochastic covariates, where $f_i$ is a known density depending on the given predictors ${\boldsymbol{x}_i}$,
error distribution, and a common regression parameter $\boldsymbol\beta$, 
$y_i\sim f_i(\cdot, {\boldsymbol{x}_i}, \boldsymbol\beta)$. 
This differs from the usual regression set-up with stochastic covariates 
that has been explored in greater detail in the literature
(\cite{Ronchetti/Roussew:1980,Schrader/Hettmansperger:1980,Ronchetti:1982a,
Ronchetti:1982b,Ronchetti:1987,Sen:1982,
Markatou/Hettmansperger:1990,Markatou/He:1994,Markatou/Manos:1996,
Cantoni/Ronchetti:2001,Liu/etc:2005,Maronna/etc:2006,Wang/Qu:2007,Hosseinian:2009,Salibian-Barrera/etc:2014}). 
Our set-up treats the regression problem from a design point of view 
where we pre-fix the covariate levels.
The robustness literature under this general I-NH set-up is limited;
some scattered attempts have been made in particular cases 
like normal regression (\cite{Huber:1983,Muller:1998}).

\cite{Ghosh/Basu:2013} proposed a global approach for estimating ${\boldsymbol{\theta}}$ 
under the I-NH set-up by minimizing the average density power divergence (DPD) measure 
(originally introduced by \cite{Basu/etc:1998}  for i.i.d. data) between the data and the model density; 
the proposed minimum DPD estimator (MDPDE) has 
excellent efficiency and robustness properties in the normal regression model. 
The approach is also implemented in the context of generalized linear models by \cite{Ghosh/Basu:2014}; 
it provides a competitive alternative to existing robust methods.
The approach has been used in \cite{Ghosh:2014} to obtain a robust alternative for 
tail index estimation under suitable assumptions of an exponential regression model.
Here, we exploit the properties of the \cite{Ghosh/Basu:2013} estimator 
to develop a general class of robust tests of hypotheses for I-NH data.


The specific advantages of the proposed methods are as follows. 
(1) The method is completely general in that it works for any set-up involving independent and non-homogeneous data. 
(2) The proposal is simple to implement with minimal addition in computational complexity
	compared to likelihood based methods. 
(3) The testing procedure is based on the minimization of a bona-fide objective function
	and the selection of the proper root of the estimating equation is simple as it must correspond 
	to the global minimum.
(4) Our methods have bounded influence for the test statistics,
	and the level and power influence functions. 
(5) The proposed tests are consistent at any fixed alternative, and 
	have high power against any contiguous alternative.

%
%
%
%

%
%

In this paper, we assume  Conditions (A1)--(A7) of \cite{Ghosh/Basu:2013},
which we refer to as the ``Ghosh-Basu conditions", and Assumptions A, B, C and D of 
\citet[][p.~429]{Lehmann:1983}, which we refer to as the ``Lehmann conditions". 
These conditions are listed in Section S1 of the Online Supplement for completeness.

\section{The MDPDE under the I-NH Set-up}
\label{SEC:8MDPDE_nonH}

Under the I-NH set-up, \cite{Ghosh/Basu:2013} proposed the estimation of ${\boldsymbol{\theta}}$ by minimizing 
the average DPD measure between the data and the model, or equivalently
\begin{eqnarray}
H_n({\boldsymbol{\theta}}) = \frac{1}{n} \sum_{i=1}^n \left[ \int f_i(y;{\boldsymbol{\theta}})^{1+\tau} dy 
- \left(1+\frac{1}{\tau}\right) f_i(Y_i;{\boldsymbol{\theta}})^\tau \right] 
= \frac{1}{n} \sum_{i=1}^n V_i(Y_i;{\boldsymbol{\theta}}).\label{EQ:8Hn}
\end{eqnarray}
The corresponding estimating equation is given by
\begin{equation}
\sum_{i=1}^n \left[ f_i(Y_i;{\boldsymbol{\theta}})^\tau \boldsymbol{u}_i(Y_i;{\boldsymbol{\theta}}) 
- \int f_i(y;{\boldsymbol{\theta}})^{1+\tau} \boldsymbol{u}_i(y;{\boldsymbol{\theta}}) dy \right] = 0, \label{EQ:general-equation}
\end{equation}
where $\nabla$ represents the gradient with respect to ${\boldsymbol{\theta}}$, and 
$\boldsymbol{u}_i(y; {\boldsymbol{\theta}}) = \nabla \ln f_i(y; {\boldsymbol{\theta}})$ 
is the likelihood score function for $i$-th model density 
(Similarly, $\nabla^2$ represents the second order derivative with respect to ${\boldsymbol{\theta}}$).
When  $\tau = 0$, the MDPDE is 
seen to coincide with the non-robust maximum likelihood estimator (MLE); 
as $\tau$ increases the robustness increases significantly at the cost of a slight loss in asymptotic efficiency.

With $\underline{\boldsymbol{G}}= (G_1,\cdots,G_n)$, 
the minimum DPD functional ${\boldsymbol{\theta}}^g = \boldsymbol{U}_\tau(\underline{\boldsymbol{G}})$ for 
the I-NH observations is defined by 
\begin{equation}
\frac{1}{n} \sum_{i=1}^n d_\tau(g_i(.),f_i(.;\boldsymbol{U}_\tau(\underline{\boldsymbol{G}}))) 
= \min_{{\boldsymbol{\theta}} \in \Theta} \frac{1}{n} \sum_{i=1}^n d_\tau(g_i(.),f_i(.;{\boldsymbol{\theta}})),
\label{EQ:MDPDE_Func}
\end{equation}
where $d_\tau(f_1, f_2)$ denotes the DPD measure between two densities $f_1$ and $f_2$ 
with the tuning parameter $\tau$, as given by \cite{Basu/etc:1998}, 
\begin{equation}\label{EQ:dpd}
d_\tau(f_1,f_2) = \displaystyle \left\{\begin{array}{ll}
\displaystyle \int  \left[f_2^{1+\tau} - \left(1 + \frac{1}{\tau}\right)  f_2^\tau f_1 + 
\frac{1}{\tau} f_1^{1+\tau}\right], & {\rm for} ~\tau > 0,\\
\displaystyle \int f_1 \log(f_1/f_2), & {\rm for} ~\tau = 0.  
\end{array}\right.
\end{equation}
Equivalently, $\boldsymbol{U}_\tau(\underline{\boldsymbol{G}})$ 
is the minimizer of $\frac{1}{n} \sum_{i=1}^n H^{(i)}({\boldsymbol{\theta}})$,
with respect to ${\boldsymbol{\theta}}\in\Theta$, where
$H^{(i)}({\boldsymbol{\theta}})= \int f_i(y;{\boldsymbol{\theta}})^{1+\tau} dy - 
\left(1+\frac{1}{\tau}\right) \int f_i(y;{\boldsymbol{\theta}})^\tau g_i(y) dy.$

\cite{Ghosh/Basu:2013} derived the asymptotic distribution of the MDPDE $\hat{{\boldsymbol{\theta}}}_n$, 
under this set-up.	Under the Ghosh-Basu conditions,
we have the following. \\
\noindent(i) There exists a consistent sequence $\widehat{{\boldsymbol{\theta}}}_n$ of roots 
of (\ref{EQ:general-equation}).
	
\noindent (ii) The asymptotic distribution of 
	$\boldsymbol\Omega_n^\tau({\boldsymbol{\theta}}^g)^{-\frac{1}{2}}\boldsymbol\Psi_n^\tau({\boldsymbol{\theta}}^g) 
	[\sqrt n (\widehat{\boldsymbol{\theta}}_n - {\boldsymbol{\theta}}^g)]$ is $p$-dimensional normal 
	with (vector) mean $\boldsymbol{0}$ and covariance matrix $\boldsymbol{I}_p$, 
	the $p$-dimensional identity matrix, where 
$\boldsymbol\Psi_n^\tau({\boldsymbol{\theta}}^g) = \frac{1}{n} \sum_{i=1}^n \boldsymbol{J}^{(i)}({\boldsymbol{\theta}}^g),$
	with	
	\begin{align}
	\boldsymbol{J}^{(i)}({\boldsymbol{\theta}}^g) 
	=& \int {\boldsymbol{u}_{i}}(y;{\boldsymbol{\theta}}^g) {\boldsymbol{u}_{i}}^T(y;{\boldsymbol{\theta}}^g) f_i^{1+\tau}(y;{\boldsymbol{\theta}}^g) dy \notag  \\
	& - \int \{ \mathbf{\nabla} {\boldsymbol{u}_{i}}(y;{\boldsymbol{\theta}}^g)  + \tau  {\boldsymbol{u}_{i}}(y;{\boldsymbol{\theta}}^g) {\boldsymbol{u}_{i}}^T(y;{\boldsymbol{\theta}}^g)\}\{g_i(y) - f_{i}(y;{\boldsymbol{\theta}}^g)\} f_{i}(y;{\boldsymbol{\theta}}^g)^\tau  dy,\notag
	\end{align}
	\begin{align}
\mbox{and}~~~
	\boldsymbol\Omega_n^\tau({\boldsymbol{\theta}}^g) 
	=& \frac{1}{n} \sum_{i=1}^n \left[ \int {\boldsymbol{u}_{i}}(y;{\boldsymbol{\theta}}^g) {\boldsymbol{u}_{i}}^T(y;{\boldsymbol{\theta}}^g) 
	f_{i}(y;{\boldsymbol{\theta}}^g)^{2\tau} g_i(y)dy 
	- \boldsymbol\xi_i  \boldsymbol\xi_i^{T} \right],\nonumber
	\end{align}
	\begin{align}
\mbox{with}~~~~~~~~~~~~~~~~~~~~~~~~	
\boldsymbol\xi_i = \int {\boldsymbol{u}_{i}}(y;{\boldsymbol{\theta}}^g) f_{i}(y;{\boldsymbol{\theta}}^g)^\tau g_i(y)dy.
~~~~~~~~~~~~~~
	\label{EQ:xi}
	\end{align}


\section{Testing Simple Hypothesis}
\label{SEC:8simple_testing}

We start with the simple hypothesis testing problem with a fully specified null under the I-NH set-up. 
Let ${\boldsymbol{\theta}}_0$ be a fixed point in the parameter space $\Theta$.  
We want to test 
\begin{equation}
H_0 : {\boldsymbol{\theta}} = {\boldsymbol{\theta}}_0 ~~~\mbox{ against }~~~~ H_1 : {\boldsymbol{\theta}} \ne {\boldsymbol{\theta}}_0.
\label{EQ:simple_Hyp}
\end{equation}
When the model is correctly specified and the null hypothesis is correct, 
$f_i(\cdot;{\boldsymbol{\theta}}_0)$ is the data generating density for the $i$-th observation. 
We can test for the hypothesis in (\ref{EQ:simple_Hyp}) by using the 
DPD measure between $f_i(\cdot;{\boldsymbol{\theta}}_0)$ and $f_i(\cdot;\widehat{\boldsymbol{\theta}})$ 
for any estimator $\widehat{{\boldsymbol{\theta}}}$ of ${\boldsymbol{\theta}}$. We consider the 
MDPDE ${\boldsymbol{\theta}}_n^\tau$ of ${\boldsymbol{\theta}}$ as defined in Section \ref{SEC:8MDPDE_nonH}.
Since there are $n$ divergence measures corresponding to each $i$, we 
consider the total divergence measure over the $n$  data points for testing (\ref{EQ:simple_Hyp})
and define the DPD based test statistic (DPDTS) as 
\begin{equation}
T_{\gamma}( {{\boldsymbol{\theta}}_n^\tau}, {{\boldsymbol{\theta}}_0}) 
= 2 \sum_{i=1}^n~d_\gamma(f_i(.;{\boldsymbol{\theta}}_n^\tau),f_i(.;{\boldsymbol{\theta}}_0)),\nonumber
\end{equation}
where $d_\gamma(f_1, f_2)$ 
is defined in (\ref{EQ:dpd}).
In case of i.i.d.~data, this DPDTS coincides with the test statistic in \cite{Basu/etc:2013a}.

\subsection{Asymptotic Properties}\label{SEC:8asymp_simple_test}

Consider the matrices $\boldsymbol\Psi_n^\tau$ and $\boldsymbol\Omega_n^\tau$ defined in 
Section \ref{SEC:8MDPDE_nonH}
and let $\boldsymbol{A}_n^\gamma({\boldsymbol{\theta}}) = \frac{1}{n} \sum_{i=1}^n \boldsymbol{A}_\gamma^{(i)}({\boldsymbol{\theta}})$
with  
$
\boldsymbol{A}_\gamma^{(i)}({\boldsymbol{\theta}}_0) = \nabla^2 d_\gamma(f_i(.;{\boldsymbol{\theta}}),f_i(.;{\boldsymbol{\theta}}_0))\big|_{{\boldsymbol{\theta}} = {\boldsymbol{\theta}}_0}. 
$
For some $p \times p$ matrices ${\boldsymbol{J}_{\tau}}$, ${\boldsymbol{V}_{\tau}}$, $A_\tau$, and ${\boldsymbol{\theta}}\in\Theta$, 
consider the assumptions.
\begin{itemize}
\item[(C1)]$\boldsymbol\Psi_n^\tau({\boldsymbol{\theta}}) \rightarrow {\boldsymbol{J}_{\tau}}({\boldsymbol{\theta}})$ and
$\boldsymbol\Omega_n^\tau({\boldsymbol{\theta}}) \rightarrow {\boldsymbol{V}_{\tau}}({\boldsymbol{\theta}})$ element-wise as $n \rightarrow \infty$
\item[(C2)]$\boldsymbol{A}_n^\gamma({\boldsymbol{\theta}}) \rightarrow {\boldsymbol{A}_{\gamma}}({\boldsymbol{\theta}})$ element-wise as $n \rightarrow \infty$. 
\end{itemize}

\begin{theorem}
Suppose the model density satisfies the Lehmann and Ghosh-Basu conditions
and conditions (C1) and (C2) hold at ${\boldsymbol{\theta}}={\boldsymbol{\theta}}_0$. 
Then, the asymptotic null distribution of the DPDTS 
$T_{\gamma}( {{\boldsymbol{\theta}}_n^\tau}, {{\boldsymbol{\theta}}_0})$ is the distribution of 
$\sum_{i=1}^r ~  \zeta_i^{\gamma, \tau}({\boldsymbol{\theta}}_0)Z_i^2,$
where $Z_1, \cdots,Z_r$ are independent standard normal variables and
$\zeta_1^{\gamma, \tau}({\boldsymbol{\theta}}_0), \cdots, \zeta_r^{\gamma, \tau}({\boldsymbol{\theta}}_0)$ are the nonzero eigenvalues
 of ${\boldsymbol{A}_{\gamma}}({\boldsymbol{\theta}}_0){\boldsymbol{\Sigma}_{\tau}}({\boldsymbol{\theta}}_0)$ with 
 ${\boldsymbol{\Sigma}_{\tau}}({\boldsymbol{\theta}}) = {\boldsymbol{J}_{\tau}}^{-1}({\boldsymbol{\theta}}) {\boldsymbol{V}_{\tau}}({\boldsymbol{\theta}}) {\boldsymbol{J}_{\tau}}^{-1}({\boldsymbol{\theta}})$ and
$r = rank( {\boldsymbol{V}_{\tau}}({\boldsymbol{\theta}}_0) {\boldsymbol{J}_{\tau}}^{-1}({\boldsymbol{\theta}}_0){\boldsymbol{A}_{\gamma}}({\boldsymbol{\theta}}_0) {\boldsymbol{J}_{\tau}}^{-1}({\boldsymbol{\theta}}_0)
 {\boldsymbol{V}_{\tau}}({\boldsymbol{\theta}}_0)).$
\label{THM:8asymp_null_simple_test}
\end{theorem}

The null distribution of the proposed DPDTS has the same form as 
that in \cite{Basu/etc:2013a,Basu/etc:2013b} for i.i.d.~observations. 
The critical region of our proposal can be easily determined from 
the relevant discussion in \cite{Basu/etc:2013a,Basu/etc:2013b}.


Next we present an approximation to its power function. 
Let $\boldsymbol{M}_n^\gamma({\boldsymbol{\theta}}) = n^{-1} \sum_{i=1}^n \boldsymbol{M}_\gamma^{(i)}({\boldsymbol{\theta}})$
with $\boldsymbol{M}_\gamma^{(i)}({\boldsymbol{\theta}})$ $= \nabla d_\gamma(f_i(.;{\boldsymbol{\theta}}),f_i(.;{\boldsymbol{\theta}}_0))$,
and assume 
\begin{itemize}
\item[(C3)]$\boldsymbol{M}_n^\gamma({\boldsymbol{\theta}}) \rightarrow 
\boldsymbol{M}_\gamma({\boldsymbol{\theta}})$ element-wise as $n \rightarrow \infty$ 
for some $p$-vector  $\boldsymbol{M}_\gamma({\boldsymbol{\theta}})$.
\end{itemize}

\begin{theorem}
Suppose the model density satisfies the Lehmann and Ghosh-Basu conditions and take any
${\boldsymbol{\theta}}^* \ne {\boldsymbol{\theta}}_0$ in $\Theta$ for which (C1) and (C3) hold.
Then, an approximation to the power function of the test 
$\left\{T_{\gamma}( {{\boldsymbol{\theta}}_n^\tau}, {{\boldsymbol{\theta}}_0}) > t_\alpha^{\tau,\gamma} \right\} $ 
for testing the hypothesis in (\ref{EQ:simple_Hyp}) at the significance level $\alpha$ is given by
\begin{eqnarray}
\pi_{n,\alpha}^{\tau,\gamma} ({\boldsymbol{\theta}}^*) = 1 - \Phi \left( \frac{1}{\sqrt{n} \sigma_{\tau,\gamma}({\boldsymbol{\theta}}^*)} 
\left( \frac{t_\alpha^{\tau,\gamma}}{2} -  \sum_{i=1}^n d_\gamma(f_i(.;{\boldsymbol{\theta}}^*),f_i(.;{\boldsymbol{\theta}}_0))\right) 
\right),\nonumber
\end{eqnarray}
where $t_\alpha^{\tau,\gamma}$ is the $(1-\alpha)$-th quantile of the asymptotic null distribution of
$T_{\gamma}({{\boldsymbol{\theta}}_n^\tau}, {{\boldsymbol{\theta}}_0})$ and 
$\sigma_{\tau,\gamma}^2({\boldsymbol{\theta}}) = \boldsymbol{M}_{\gamma}({\boldsymbol{\theta}})^T {\boldsymbol{\Sigma}_{\tau}}({\boldsymbol{\theta}})\boldsymbol{M}_{\gamma}({\boldsymbol{\theta}}).$
\label{THM:8asymp_power_simple_test}
\end{theorem}
\begin{corollary}
	For any ${\boldsymbol{\theta}}^* \ne {\boldsymbol{\theta}}_0$, the probability of rejecting the null hypothesis $H_0$ at any fixed 
	significance level $\alpha > 0 $ with the rejection rule $\left\{T_{\gamma}( {{\boldsymbol{\theta}}_n^\tau}, {{\boldsymbol{\theta}}_0}) >
	t_\alpha^{\tau,\gamma}\right\}$ tends to $1$ as $n \rightarrow \infty$, provided 
	$\frac{1}{n} \sum_{i=1}^n d_\gamma(f_i(.;{\boldsymbol{\theta}}^*),f_i(.;{\boldsymbol{\theta}}_0))$ $= O(1) $. 
\end{corollary}

	Theorem \ref{THM:8asymp_power_simple_test} can be used to obtain 
	the sample size required to achieve a pre-specified power $\eta$. 
	For this we just need to solve the equation
	\begin{eqnarray}
	\eta = 1 - \Phi \left( \frac{1}{\sqrt{n} \sigma_{\tau,\gamma}({\boldsymbol{\theta}}^*)} 
	\left(\frac{t_\alpha^{\tau,\gamma}}{2}-\sum_{i=1}^n d_\gamma(f_i(.;{\boldsymbol{\theta}}^*),f_i(.;{\boldsymbol{\theta}}_0))\right)\right)\nonumber
	\end{eqnarray}
	in terms of $n$. 
	If $n^*$ denotes the solution, then the required sample size is 
	the least integer greater than or equal to $n^*$.

\subsection{Robustness Properties}\label{SEC:8Robust_simple_test}

\subsubsection{Influence Functions of the Test Statistics}\label{SEC:8IF_simple_test}

We illustrate the robustness of the proposed DPDTS using some extensions of  
Hampel's influence function (IF), 
as in \cite{Huber:1983} and \cite{Ghosh/Basu:2013}.

Ignoring the multiplier $2$ in the DPDTS, we consider the functional 
$$
T_{\gamma, \tau}^{(1)}(\underline{\mathbf{G}}) 
= \sum_{i=1}^n ~ d_\gamma(f_i(\cdot;\boldsymbol{U}_\tau(\underline{\mathbf{G}})),f_i(\cdot;{\boldsymbol{\theta}}_0)), $$
where $\underline{\mathbf{G}}=(G_1,\cdots,G_n)$ and $\boldsymbol{U}_\tau(\underline{\mathbf{G}})$ is 
the minimum DPD functional under the I-NH set-up.
Unlike the i.i.d.~case, here the functional depends on the sample size $n$,
and so does the corresponding IF. 
We refer to it as the fixed-sample IF. 
Consider the contaminated distribution  $G_{i,\epsilon} = (1-\epsilon) G_i + \epsilon \wedge_{t_i}$, 
where $\wedge_{t_i}$ is the degenerate distribution at the point of contamination $t_i$ in the 
$i$-th direction for all $i=1, \ldots, n$. 
Thus, we can have contamination in some fixed direction or in all directions.

Consider a contamination only in the $i_0$-th direction and take 
$\underline{\mathbf{G}}_{i_0,\epsilon} = (G_1$, $\cdots$, $ G_{i_0-1}, G_{i_0,\epsilon}$, $\cdots, G_n) $. 
Then the corresponding first order IF of the test functional 
$T_{\gamma, \tau}^{(1)}(\underline{\mathbf{G}})$ can be defined as 
\begin{eqnarray}
IF_{i_0}(t_{i_0}, T_{\gamma, \tau}^{(1)}, \underline{\mathbf{G}})  &=& \left.\frac{\partial}{\partial\epsilon} 
T_{\gamma, \tau}^{(1)}(\underline{\mathbf{G}}_{i_0,\epsilon}) \right|_{\epsilon=0} 
=\sum_{i=1}^{n} \boldsymbol{M}_\gamma^{(i)}(\boldsymbol{U}_\tau(\underline{\mathbf{G}}))^T 
	IF_{i_0}(t_{i_0}, \boldsymbol{U}_\tau, \underline{\mathbf{G}}),\nonumber
\end{eqnarray}
where $ IF_{i_0}(t_{i_0}, \boldsymbol{U}_\tau, \underline{\mathbf{G}})$ is the IF of 
$\boldsymbol{U}_{\tau}$ in \cite{Ghosh/Basu:2013}, 
\begin{align}
IF_{i_0}(t_{i_0}, \boldsymbol{U}_\tau, \underline{\mathbf{G}})
=& \frac{1}{n} \boldsymbol\Psi_n^\tau({\boldsymbol{\theta}^g})^{-1}\boldsymbol{D}_{\tau, i_0}(t_{i_0};{\boldsymbol{\theta}^g}),
\label{EQ:IF_MDPDE}
\end{align}
where 
$\boldsymbol{D}_{\tau, i}(t;{\boldsymbol{\theta}})
=\left[  f_{i}(t;{\boldsymbol{\theta}})^\tau {\boldsymbol{u}_{i}}(t;{\boldsymbol{\theta}}) - \boldsymbol\xi_{i}\right]$
with $\boldsymbol\xi_{i}$ defined at (\ref{EQ:xi}).
In general, the IF of a test is evaluated at the null distribution 
$G_i(\cdot)=F_i(\cdot,{{\boldsymbol{\theta}}_0})$ for all $i$. 
Letting $\mathbf{\underline{F}}_{{\boldsymbol{\theta}}_0}$$=(F_1(\cdot,{{\boldsymbol{\theta}}_0}),$$ \cdots, $$F_n(\cdot,{{\boldsymbol{\theta}}_0})),$ 
we get	$\boldsymbol{U}_\tau(\mathbf{\underline{F}}_{{\boldsymbol{\theta}}_0}) = {\boldsymbol{\theta}}_0$ and 
$\boldsymbol{M}_\gamma^{(i)}({\boldsymbol{\theta}}_0)=\boldsymbol{0}$ so that  Hampel's first-order IF of  
the DPDTS is zero at $H_0$.

The second order IF of the DPDTS can be defined similarly as 
\begin{eqnarray}
IF_{i_0}^{(2)}(t_{i_0}, T_{\gamma, \tau}^{(1)}, \underline{\mathbf{G}}) 
= \frac{\partial^2}{\partial^2\epsilon} 
T_{\gamma, \tau}^{(1)}(G_1,\cdots,G_{i_0-1},G_{i_0,\epsilon},\cdots, G_n) \big|_{\epsilon=0}. 
\nonumber 
\end{eqnarray}
In particular, at the null distribution  $\underline{\mathbf{G}}=\mathbf{\underline{F}}_{{\boldsymbol{\theta}}_0}$, 
it simplifies to
\begin{equation}
IF_{i_0}^{(2)}(t_{i_0}, T_{\gamma, \tau}^{(1)}, \mathbf{\underline{F}}_{{\boldsymbol{\theta}}_0}) 
= n \cdot IF_{i_0}(t_{i_0}, \boldsymbol{U}_\tau, \mathbf{\underline{F}}_{{\boldsymbol{\theta}}_0})^T 
\boldsymbol{A}_n^\gamma({\boldsymbol{\theta}_0}) 
IF_{i_0}(t_{i_0}, \boldsymbol{U}_\tau, \mathbf{\underline{F}}_{{\boldsymbol{\theta}}_0}). \nonumber
\end{equation}
Thus the IF of the test at the null is bounded for any fixed sample size if and only if 
the IF of the corresponding minimum DPD functional is bounded. 
Using the form of the IF of the MDPDE from (\ref{EQ:MDPDE_Func}), 
the IF of the test is 
\begin{eqnarray}
IF_{i_0}^{(2)}(t_{i_0}, T_{\gamma, \tau}^{(1)}, \mathbf{\underline{F}}_{{\boldsymbol{\theta}}_0}) 
&=& \frac{1}{n} \boldsymbol{D}_{\tau, i_0}(t_{i_0};{\boldsymbol{\theta}}_0)^T 
[\boldsymbol\Psi_n^\tau({\boldsymbol{\theta}_0})^{-1} \boldsymbol{A}_n^\gamma({\boldsymbol{\theta}_0}) 
\boldsymbol\Psi_n^\tau({\boldsymbol{\theta}_0})^{-1}] 
\boldsymbol{D}_{\tau, i_0}(t_{i_0};{\boldsymbol{\theta}}_0).\nonumber
\end{eqnarray}
Here $\boldsymbol{D}_{\tau, i}(t;{\boldsymbol{\theta}})$ is bounded in $t$ if the parametric model satisfies  
$f_{i}(t;{\boldsymbol{\theta}})^\tau {\boldsymbol{u}_{i}}(t;{\boldsymbol{\theta}})$ is bounded in $t$.  
For most parametric models, this holds at $\tau>0$ (but not at $\tau=0$)
implying that the $\boldsymbol{D}_{\tau, i}(t;{\boldsymbol{\theta}})$, 
and therefore the IF, is bounded whenever $\tau>0$, but unbounded at $\tau=0$.
Note that $\boldsymbol{D}_{\tau, i}(t;{\boldsymbol{\theta}})$ does not depend on the tuning parameter $\gamma$,
hence its boundedness and that of the IF of the proposed test 
is independent of the choice of $\gamma$.

If we consider contamination in all directions at the contamination 
point $\mathbf{t}=(t_1, \cdots, t_n)$, we can derive 
the corresponding IF of the proposed DPDTS in a similar manner. 
Again, at the null distribution, its first order IF is zero 
and its second order IF simplifies to 
\begin{eqnarray}
&& IF^{(2)}(\mathbf{t}, T_{\gamma, \tau}^{(1)}, \mathbf{\underline{F}}_{{\boldsymbol{\theta}}_0}) 
= n \cdot IF(\mathbf{t}, \boldsymbol{U}_\tau, \mathbf{\underline{F}}_{{\boldsymbol{\theta}}_0})^T A_n^\gamma 
IF(\mathbf{t}, \boldsymbol{U}_\tau, \mathbf{\underline{F}}_{{\boldsymbol{\theta}}_0}).\nonumber\\
&& ~~= \frac{1}{n} \left( \sum_{i=1}^n \boldsymbol{D}_{\tau, i}(t_i;{\boldsymbol{\theta}}_0)\right)^T 
[\boldsymbol\Psi_n^\tau({\boldsymbol{\theta}_0})^{-1} \boldsymbol{A}_n^\gamma({\boldsymbol{\theta}_0}) 
\boldsymbol\Psi_n^\tau({\boldsymbol{\theta}_0})^{-1} ] 
\left( \sum_{i=1}^n \boldsymbol{D}_{\tau, i}(t_i;{\boldsymbol{\theta}}_0)\right).\nonumber
\end{eqnarray}
This influence function is bounded for most parametric models when $\tau>0$ 
and unbounded if $\tau=0$. Thus, whatever be the contamination direction, 
the proposed DPDTS is always robust for $\tau>0$ and non-robust for $\tau=0$.
Here, robustness refers to local robustness of the test statistics
under infinitesimal contamination.

\subsubsection{Level and Power Influence Functions}
\label{SEC:8IF_power_simple_test}

The performance of any testing procedure is generally measured by its level and power.
We consider the effect of contamination on level and power of the proposed DPDTS
through the level and power influence functions (\cite{Hampel/etc:1986,Heritier/Ronchetti:1994,Toma/Broniatowski:2010}).
Since the exact level and power of the proposed test are difficult to obtain, we work with their asymptotic versions.

Since the proposed DPDTS is consistent, we examine its asymptotic power under the contiguous alternative 
$H_{1,n} : {\boldsymbol{\theta}}_n = {\boldsymbol{\theta}}_0 + n^{-1/2}\boldsymbol\Delta$ with $\boldsymbol\Delta \in \mathbb{R}^p - \{0\}$. 
Here we consider contamination over these alternatives. 
	As argued in \cite{Hampel/etc:1986}, we consider contaminations such that their effects tend to zero as 
	${\boldsymbol{\theta}}_n$ tends to ${\boldsymbol{\theta}}_0$ at the same rate to avoid the confusion between the null and alternative 
	neighborhoods. 
Consider the contaminated distributions 
	$$
	\mathbf{\underline{F}}_{n,\epsilon,\mathbf{t}}^L = \left(1-\frac{\epsilon}{\sqrt{n}}\right) 
	\mathbf{\underline{F}}_{{\boldsymbol{\theta}}_0}+\frac{\epsilon}{\sqrt{n}} \wedge_{\mathbf{t}},
	\mbox{ and }
	\mathbf{\underline{F}}_{n,\epsilon,\mathbf{t}}^P = \left(1-\frac{\epsilon}{\sqrt{n}}\right) 
	\mathbf{\underline{F}}_{{\boldsymbol{\theta}}_n}+\frac{\epsilon}{\sqrt{n}} \wedge_{\mathbf{t}},
	$$
for the level and power, respectively,
	where $\mathbf{t} = (t_1, \cdots, t_n)^T$, 
	$\mathbf{\underline{F}}_{n,\epsilon,\mathbf{t}}^P = (F_{i,n,\epsilon,t_i}^P)_{i=1,\cdots,n}$ and 
	$\mathbf{\underline{F}}_{n,\epsilon,\mathbf{t}}^L = (F_{i,n,\epsilon,t_i}^L)_{i=1,\cdots,n}$. 
	Then the level influence function (LIF) and the power influence function (PIF) 
	are
\begin{eqnarray}
	LIF(\mathbf{t}; T_{\gamma}^{(1)}, \mathbf{\underline{F}}_{{\boldsymbol{\theta}}_0} ) 
	&=& \lim_{n \rightarrow \infty} ~ \frac{\partial}{\partial \epsilon} 
	P_{\mathbf{\underline{F}}_{n,\epsilon,\mathbf{t}}^L }( T_{\gamma}( {{\boldsymbol{\theta}}_n^\tau}, {{\boldsymbol{\theta}}_0})
	>  t_\alpha^{\tau,\gamma}) \big|_{\epsilon=0},\nonumber\\
	PIF(\mathbf{t}; T_{\gamma}^{(1)}, \mathbf{\underline{F}}_{{\boldsymbol{\theta}}_0} ) 
	&=& \lim_{n \rightarrow \infty} ~ \frac{\partial}{\partial \epsilon} 
	P_{\mathbf{\underline{F}}_{n,\epsilon,\mathbf{t}}^P }( T_{\gamma}( {{\boldsymbol{\theta}}_n^\tau}, {{\boldsymbol{\theta}}_0}) 
	>  t_\alpha^{\tau,\gamma}) \big|_{\epsilon=0}.\nonumber
\end{eqnarray}
We start with the asymptotic power under the contaminated distribution 
$\mathbf{\underline{F}}_{n,\epsilon,\mathbf{y}}^P$ and examine some special cases
by substituting specific values of $\boldsymbol\Delta$ and $\epsilon$.

\begin{theorem}
	Suppose that the Lehmann and Ghosh-Basu conditions hold for the model density
	and (C1)-(C2) hold at ${\boldsymbol{\theta}}={\boldsymbol{\theta}}_0$. 
	Then for any $\boldsymbol\Delta \in \mathbb{R}^p$ and $\epsilon \geq 0$, we have the following.

\noindent(i) The asymptotic  distribution of the proposed DPDTS
		under $\mathbf{\underline{F}}_{n,\epsilon,\mathbf{t}}^P $ 
		is the distribution of the quadratic form 
		$\boldsymbol{W}^T {\boldsymbol{A}_{\gamma}}({\boldsymbol{\theta}}_0)\boldsymbol{W}$, 
		where $\boldsymbol{W}\sim N_p\left(\widetilde{\boldsymbol\Delta}, {\boldsymbol{\Sigma}_{\tau}}({\boldsymbol{\theta}}_0)\right)$ 
		with $\widetilde{\boldsymbol\Delta} = \left[\boldsymbol\Delta 
		+ \epsilon IF(\mathbf{t}; \boldsymbol{U}_\tau, \mathbf{\underline{F}}_{{\boldsymbol{\theta}}_0})\right]$
		and $\boldsymbol{\Sigma}_\tau$ defined as in Theorem \ref{THM:8asymp_null_simple_test}. 	
		
\noindent(ii) The asymptotic  power of the proposed DPDTS under 
		$\mathbf{\underline{F}}_{n,\epsilon,\mathbf{t}}^P $  is given by
		\begin{eqnarray}
		P_{\tau,\gamma}(\boldsymbol\Delta, \epsilon; \alpha) &=& 
		\lim_{n \rightarrow \infty} ~ P_{\mathbf{\underline{F}}_{n,\epsilon,\mathbf{t}}^L }( 
		T_{\gamma}( {{\boldsymbol{\theta}}_n^\tau}, {{\boldsymbol{\theta}}_0}) >  t_\alpha^{\tau,\gamma})
		,\nonumber\\&=& 
		 \sum\limits_{v=0}^{\infty}~C_v^{\gamma,\tau}({\boldsymbol{\theta}}_0, \widetilde{\boldsymbol\Delta})
		P\left(\chi_{r+2v}^2 > \frac{t_\alpha^{\tau, \gamma}}{\zeta_{(1)}^{\gamma,\tau}({\boldsymbol{\theta}}_0)}\right), 
		\label{EQ:8asymp_power_cont_null}
		\end{eqnarray}
		where  $\chi_p^2$ denotes  a  chi-square  random  variable with $p$  degrees of  freedom,
		$\zeta_{(1)}^{\gamma,\tau}({\boldsymbol{\theta}}_0)$ is the minimum of $\zeta_{i}^{\gamma,\tau}({\boldsymbol{\theta}}_0)$s 
		for $i=1,\ldots, r$ defined in Theorem \ref{THM:8asymp_null_simple_test}, and
		$$
		C_v^{\gamma,\tau}({\boldsymbol{\theta}}_0, \widetilde{\boldsymbol\Delta}) = \frac{1}{v!} \left(
		\prod\limits_{j=1}^{r}\frac{\zeta_{(1)}^{\gamma,\tau}({\boldsymbol{\theta}}_0)}{
			\zeta_{j}^{\gamma,\tau}({\boldsymbol{\theta}}_0)}\right)^{1/2}
		e^{-\frac{1}{2}\sum\limits_{j=1}^{r}\delta_j} E\left[\left(\widehat{Q}\right)^v\right], \nonumber
		$$
		$$
\mbox{with }~~~~~~~
		\widehat{Q} = \frac{1}{2}\sum\limits_{j=1}^{r}\left[\left(1 - 
		\frac{\zeta_{(1)}^{\gamma,\tau}({\boldsymbol{\theta}}_0)}{\zeta_{j}^{\gamma,\tau}({\boldsymbol{\theta}}_0)}\right)^{1/2}Z_j
		+ \sqrt{\delta_j}\left(\frac{\zeta_{(1)}^{\gamma,\tau}({\boldsymbol{\theta}}_0)}{\zeta_{j}^{\gamma,\tau}({\boldsymbol{\theta}}_0)}\right)^{1/2}
		\right]^2,
		$$
		for $r$ independent standard normal random variables $Z_1, \ldots, Z_r$,
		and $\delta_i$s defined as in Remark \ref{REM:a} below.
	\label{THM:8asymp_Contaminated_power_simple_test}
\end{theorem}

\begin{corollary}\label{COR:8contiguous_power_simple_test}
The asymptotic power under the contiguous alternatives 
	$H_{1,n}: {\boldsymbol{\theta}}= {\boldsymbol{\theta}}_n = {\boldsymbol{\theta}}_0 + n^{-\frac{1}{2}}\boldsymbol\Delta$ is 
	\begin{eqnarray}
	P_{\tau,\gamma}(\boldsymbol\Delta, 0; \alpha) &=& 
	\sum\limits_{v=0}^{\infty}~C_v^{\gamma,\tau}({\boldsymbol{\theta}}_0, {\boldsymbol\Delta})
	P\left(\chi_{r+2v}^2 > \frac{t_\alpha^{\tau, \gamma}}{\zeta_{(1)}^{\gamma,\tau}({\boldsymbol{\theta}}_0)}\right). \nonumber
	\end{eqnarray}
\end{corollary}

\begin{corollary}\label{COR:8asymp_level_simple_test}
The asymptotic level under the probability distribution $\mathbf{\underline{F}}_{n,\epsilon,\mathbf{t}}^L$ is 
	\begin{eqnarray}
	\alpha_\epsilon = P_{\tau,\gamma}(\boldsymbol{0}, \epsilon; \alpha) &=& 
	\sum\limits_{v=0}^{\infty}~C_v^{\gamma,\tau}\left({\boldsymbol{\theta}}_0, 
	\epsilon IF(\mathbf{t}; \boldsymbol{U}_\tau, \mathbf{\underline{F}}_{{\boldsymbol{\theta}}_0})\right)
	P\left(\chi_{r+2v}^2 > \frac{t_\alpha^{\tau, \gamma}}{\zeta_{(1)}^{\gamma,\tau}({\boldsymbol{\theta}}_0)}\right). \nonumber
	\end{eqnarray}
\end{corollary}

\begin{remark}
\label{REM:a}		
The asymptotic  distribution of $T_{\gamma}( {{\boldsymbol{\theta}}_n^\tau}, {{\boldsymbol{\theta}}_0})$ 
under $\mathbf{\underline{F}}_{n,\epsilon,\mathbf{t}}^P $
is the same as that of $\sum\limits_{i=1}^{r}\zeta_i^{\gamma, \tau}({\boldsymbol{\theta}}_0)\chi_{1,\delta_i}^2$,  where 
the	$\zeta_i^{\gamma, \tau}({\boldsymbol{\theta}}_0)$ are	given in Theorem \ref{THM:8asymp_null_simple_test} and 
		$\chi_{1,\delta_i}^2$s	are independent  non-central  chi-square  variables 
		having degree of freedom one and non-centrality	parameters $\delta_i$s, respectively, with
		$\left(\sqrt{\delta_1}, \ldots, \sqrt{\delta_p}\right)^T 
		= \widetilde{\boldsymbol{P}}_{\tau,\gamma}({\boldsymbol{\theta}}_0){\boldsymbol{\Sigma}_{\tau}}^{-1/2}
		({\boldsymbol{\theta}}_0)\widetilde{\boldsymbol\Delta}$
		and $\widetilde{\boldsymbol{P}}_{\tau,\gamma}({\boldsymbol{\theta}}_0)$ the matrix of normalized eigenvectors of 
		${\boldsymbol{A}_{\gamma}}({\boldsymbol{\theta}}_0){\boldsymbol{\Sigma}_{\tau}}({\boldsymbol{\theta}}_0)$.
\end{remark}

\begin{remark}
The expressions of asymptotic level and power 
under contiguous alternative with contamination can be approximated
by truncating the series to a finite number ($N$) terms.  
The error incurred by such a truncation can be made smaller than any pre-specific limit 
by choosing $N$ suitably large.
If we truncate at the $N$-th term of (\ref{EQ:8asymp_power_cont_null}),
assuming $c_v = C_v^{\gamma, \tau}({\boldsymbol{\theta}}_0, \widetilde{\boldsymbol\Delta})$,
the error can be bounded by 
\begin{eqnarray}
e_N &=& \sum\limits_{v=N+1}^{\infty} ~ c_v\cdot
P\left(\chi_{r+2v}^2 > \frac{t_\alpha^{\tau,\gamma}}{\zeta_{(1)}^{\gamma, \tau}({\boldsymbol{\theta}}_0)}\right) 
\leq\sum\limits_{v=N+1}^{\infty} ~ c_v = 1 - \sum\limits_{v=0}^{N} ~ c_v. \nonumber
\end{eqnarray}
See \cite{Kotz/etc:1967a,Kotz/etc:1967b} for more accurate error bounds for such approximations.
\end{remark}

Starting with $P_{\tau,\gamma}(\boldsymbol\Delta, \epsilon; \alpha)$ at (\ref{EQ:8asymp_power_cont_null})
and differentiating, we get the power influence function $PIF(\cdot)$.


\begin{theorem}\label{THM:8IF_power_simple_test}
	Assume that the Lehmann and Ghosh-Basu conditions hold for the model density
	and (C1)-(C2) hold at ${\boldsymbol{\theta}}={\boldsymbol{\theta}}_0$. 
If the influence function $IF(\mathbf{t}; \boldsymbol{U}_\tau,\mathbf{\underline{F}}_{{\boldsymbol{\theta}}_0})$ of the MDPDE is bounded,
then for any $\boldsymbol\Delta \in \mathbb{R}^p$, the power influence function of the proposed DPDTS is given by 
	%
$$
PIF(\mathbf{t}; T_{\gamma,\lambda}^{(1)}, \mathbf{\underline{F}}_{{\boldsymbol{\theta}}_0}) 
= IF(\mathbf{t}; \boldsymbol{U}_\tau, \mathbf{\underline{F}}_{{\boldsymbol{\theta}}_0})^T 
\boldsymbol{K}_{\gamma,\tau}({\boldsymbol{\theta}}_0,\boldsymbol\Delta,\alpha),
$$
where
$\boldsymbol{K}_{\gamma,\tau}({\boldsymbol{\theta}}_0,\boldsymbol\Delta,\alpha)=
\left(\sum\limits_{v=0}^{\infty}~\left[\left.\frac{\partial}{\partial \boldsymbol{d}} 
	C_v^{\gamma,\tau}({\boldsymbol{\theta}}_0, \boldsymbol{d}) \right|_{\boldsymbol{d}=\boldsymbol\Delta}\right]
	P\left(\chi_{r+2v}^2 > \frac{t_\alpha^{\tau, \gamma}}{\zeta_{(1)}^{\gamma,\tau}({\boldsymbol{\theta}}_0)}\right)\right).
$
\end{theorem}

The level influence function of the proposed DPDTS can be derived by putting 
$\boldsymbol\Delta = \boldsymbol{0}$ in this expression of the PIF, yielding
$LIF(\mathbf{t}; T_{\gamma,\lambda}^{(1)}, \mathbf{\underline{F}}_{{\boldsymbol{\theta}}_0})
= IF(\mathbf{t}; \boldsymbol{U}_\tau, \mathbf{\underline{F}}_{{\boldsymbol{\theta}}_0})^T 
\boldsymbol{K}_{\gamma,\tau}({\boldsymbol{\theta}}_0,\boldsymbol{0},\alpha)$
whenever the IF of the MDPDE used is bounded. 
Thus asymptotically the level of the DPDTS is unaffected by 
contiguous contaminations for all $\tau > 0$.

\section{Testing Composite Hypothesis}\label{SEC:8composite_testing}

In the I-NH set-up of Section \ref{SEC:intro}, 
take a fixed (proper) subspace $\Theta_0$ of $\Theta$.  
Based on the observed data, we want to test the hypothesis 
\begin{eqnarray}\label{EQ:8composite_hypo}
 H_0 : {\boldsymbol{\theta}} \in \Theta_0 ~~~\mbox{ against }~~~~ H_1 : {\boldsymbol{\theta}}  \notin \Theta_0.
\end{eqnarray}
When the model is correctly specified and $H_0$ is correct, 
$f_i(\cdot;{\boldsymbol{\theta}}_0)$ is the data generating density for the $i$-th observation, 
where ${\boldsymbol{\theta}}_0 \in \Theta_0$.
We can test this hypothesis by using the DPD measure between $f_i(\cdot;\widetilde{{\boldsymbol{\theta}}})$ 
and $f_i(\cdot;\widehat{\boldsymbol{\theta}})$ for any two estimators $\widetilde{{\boldsymbol{\theta}}}$ and 
$\widehat{{\boldsymbol{\theta}}}$ of ${\boldsymbol{\theta}}$ 
under $H_0$ and $H_0 \cup H_1$, respectively. For $\widehat{{\boldsymbol{\theta}}}$,  
we take the MDPDE ${\boldsymbol{\theta}}_n^\tau$  of ${\boldsymbol{\theta}}$ with tuning parameter $\tau$, 
and for $\widetilde{{\boldsymbol{\theta}}}$, we take the estimator $\widetilde{{\boldsymbol{\theta}}}_n^\tau$ obtained by 
minimizing the DPD with tuning parameter $\tau$  over the subspace $\Theta_0$ only; 
we refer to $\widetilde{{\boldsymbol{\theta}}}_n^\tau$ as the restricted MDPDE (RMDPDE). 
Thus, our test statistic (DPDTS$_C$) for the hypothesis (\ref{EQ:8composite_hypo}) 
 based on the DPD with parameter $\gamma$  is 
\begin{eqnarray}\label{EQ:8DPDTS2}
S_{\gamma}( {{\boldsymbol{\theta}}_n^\tau}, {\widetilde{{\boldsymbol{\theta}}}_n^\tau}) 
= 2 \sum_{i=1}^n d_\gamma(f_i(.;{\boldsymbol{\theta}}_n^\tau),f_i(.;\widetilde{{\boldsymbol{\theta}}}_n^\tau)).
\end{eqnarray}

\subsection{Properties of the RMDPDE under the I-NH Set-up}
\label{SEC:8RMDPDE}

The restricted minimum density power divergence estimators (RMDPDE)  
$\widetilde{{\boldsymbol{\theta}}}_n^\tau$  of ${\boldsymbol{\theta}}$ is the minimizer of 
the DPD objective function $H_n({\boldsymbol{\theta}})$ given at (\ref{EQ:8Hn}) 
with tuning parameter $\tau$ subject to a set of $r$ restrictions of the form
$   \boldsymbol\upsilon({\boldsymbol{\theta}}) = \boldsymbol{0}$,
where $\boldsymbol\upsilon : \mathbb{R}^p \mapsto \mathbb{R}^r$ is some vector valued function.
For the null hypothesis in (\ref{EQ:8composite_hypo}), 
such restrictions are given by the definition of the null parameter space $\Theta_0$. 
We assume that the $p \times r$ matrix 
$   \boldsymbol\Upsilon({\boldsymbol{\theta}}) = \frac{\partial \boldsymbol\upsilon({\boldsymbol{\theta}})}{\partial {\boldsymbol{\theta}}} $
exists and is continuous in ${\boldsymbol{\theta}}$ with rank $r$. 
Then, the RMDPDE has to satisfy 
  \begin{eqnarray}
\nabla H_n({\boldsymbol{\theta}}) + \boldsymbol\Upsilon({\boldsymbol{\theta}})\boldsymbol\lambda_n =0,
~~~~~~~~~~~
\boldsymbol\upsilon({\boldsymbol{\theta}}) = 0,
\label{EQ:8Est_Eqn_RMDPDE}
\end{eqnarray}
where $\boldsymbol\lambda_n$ is an $r$-vector of Lagrangian multipliers.  
Further, the restricted minimum DPD functional 
$\widetilde{{\boldsymbol{\theta}}}^g = \widetilde{\boldsymbol{U}}_\tau(\underline{\mathbf{G}})$ 
at the true distribution $\underline{\mathbf{G}}$
is the minimizer of $n^{-1}\sum_{i=1}^n d_\tau(g_i(.),f_i(.;{\boldsymbol{\theta}}))$
subject to  $\boldsymbol\upsilon({\boldsymbol{\theta}})=\boldsymbol{0}$.

\begin{theorem}\label{THM:8asymp_RMDPDE}
Assume that the Ghosh-Basu conditions are satisfied with respect to $\Theta_0$ (instead of $\Theta$). 
Then the following hold. 

\noindent(i) There exists a consistent sequence $\widetilde{{\boldsymbol{\theta}}}_n^\tau$ of roots 
of (\ref{EQ:8Est_Eqn_RMDPDE}).

\noindent(ii) Asymptotically, $\boldsymbol\Omega_n^\tau(\widetilde{{\boldsymbol{\theta}}}^g)^{-\frac{1}{2}}
\boldsymbol{P}_n^\tau(\widetilde{{\boldsymbol{\theta}}}^g)^{-1}	
[\sqrt n (\widetilde{{\boldsymbol{\theta}}}_n^\tau - \widetilde{{\boldsymbol{\theta}}}^g)]
\sim N_p\left(\boldsymbol{0}, \boldsymbol{I}_p\right)$ 
where $\boldsymbol{I}_p$ is the $p\times p$ identity matrix, 
$\boldsymbol\Upsilon_n^\ast({\boldsymbol{\theta}}) = \boldsymbol\Upsilon({\boldsymbol{\theta}})^T
[\nabla^2 H_n({\boldsymbol{\theta}})]^{-1}\boldsymbol\Upsilon({\boldsymbol{\theta}})$,
and 
\begin{eqnarray}
\boldsymbol{P}_n^\tau({\boldsymbol{\theta}}) 
= \left[\frac{\nabla^2 H_n({\boldsymbol{\theta}})}{(1+\tau)}\right]^{-1} 
\left[\boldsymbol{I}_p - \boldsymbol\Upsilon({\boldsymbol{\theta}}) 
\left[\boldsymbol\Upsilon_n^\ast({\boldsymbol{\theta}})\right]^{-1} \boldsymbol\Upsilon({\boldsymbol{\theta}})^T 
[\nabla^2 H_n({\boldsymbol{\theta}})]^{-1} \right].
\nonumber\label{EQ:8P_n}
	\end{eqnarray} 
\end{theorem}

In the following, we need a further assumption. 
\begin{itemize}
\item[(C4)] $\boldsymbol{P}_n^\tau(\widetilde{{\boldsymbol{\theta}}}^g) \rightarrow \boldsymbol{P}_\tau(\widetilde{{\boldsymbol{\theta}}}^g)$ 
($p \times p$ invertible) element-wise as $n \rightarrow \infty$. 
\end{itemize}

\begin{corollary}
If the assumptions of Theorem \ref{THM:8asymp_RMDPDE} hold, and 
(C1) and (C4)  hold at ${\boldsymbol{\theta}} = \widetilde{{\boldsymbol{\theta}}}^g$,
then, asymptotically, $\sqrt n (\widetilde{{\boldsymbol{\theta}}}_n^\tau - \widetilde{{\boldsymbol{\theta}}}^g)
\sim N_p\left(\boldsymbol{0}, {\boldsymbol{P}_{\tau}}(\widetilde{{\boldsymbol{\theta}}}^g) {\boldsymbol{V}_{\tau}}(\widetilde{{\boldsymbol{\theta}}}^g) {\boldsymbol{P}_{\tau}}(\widetilde{{\boldsymbol{\theta}}}^g)\right)$. 
\end{corollary}

Next, we explore the robustness properties of the RMDPDEs in terms of their influence function. 
%
First consider the contamination in only one ($i_0$-th) direction.
Suppose the given restrictions are such that they can be substituted explicitly in the 
expression of average DPD before taking its derivative with respect to ${\boldsymbol{\theta}}$;
then the final derivative should be zero at 
${\boldsymbol{\theta}}=\widetilde{\boldsymbol{U}}_\tau(\underline{\mathbf{G}}_{i_0,\epsilon})$ 
and $g_{i_0} = g_{i_0,\epsilon}$, the density corresponding to $G_{i_0,\epsilon}$. 
Standard differentiation of the resulting equation with respect to $\epsilon$ at $\epsilon=0$ 
yields the IF of the RMDPDE,  
 $IF_{i_0}(t_{i_0}; \widetilde{\boldsymbol{U}}_\tau; \underline{\mathbf{G}})
 =\frac{\partial}{\partial\epsilon}\widetilde{\boldsymbol{U}}_\tau(\underline{\mathbf{G}}_{i_0,\epsilon})\big|_{\epsilon=0}$ 
 as a solution of  
 \begin{eqnarray}
  && \boldsymbol\Psi_n^{(0)}(\widetilde{{\boldsymbol{\theta}}^g }) 
  IF_{i_0}(t_{i_0}, \widetilde{\boldsymbol{U}}_\tau, \underline{\mathbf{G}}) - 
 \frac{1}{n} \boldsymbol{D}_{\tau, i_0}^{(0)}(t_{i_0}; \widetilde{{\boldsymbol{\theta}}^g }) = 0, 
  \label{EQ:8IF_eq1}
 \end{eqnarray}
where $\boldsymbol{D}_{\tau, i}^{(0)}(t; {\boldsymbol{\theta}}) 
= \left[f_{i}(t; {\boldsymbol{\theta}})^\tau {\boldsymbol{u}_{i}}^{(0)}(t;{\boldsymbol{\theta}}) - \boldsymbol\xi_{i}^{(0)}({\boldsymbol{\theta}})\right]$ 
and $\boldsymbol\Psi_n^{(0)}({\boldsymbol{\theta}})$, $\boldsymbol\xi_i^{(0)}({\boldsymbol{\theta}})$, 
$\boldsymbol{u}_i^{(0)}(y;{\boldsymbol{\theta}})$ are the same as 
$\boldsymbol\Psi_n({\boldsymbol{\theta}})$, $\boldsymbol\xi_i({\boldsymbol{\theta}})$, 
$\boldsymbol{u}_i(y;{\boldsymbol{\theta}})$ respectively, under the additional restriction 
$\boldsymbol\upsilon({\boldsymbol{\theta}})=0$. 
Also, $\widetilde{\boldsymbol{U}}_\tau(\underline{\mathbf{G}}_{i_0,\epsilon})$ must satisfy 
$ \boldsymbol\upsilon(\widetilde{\boldsymbol{U}}_\tau(\underline{\mathbf{G}}_{i_0,\epsilon})) = \boldsymbol{0}$,
from which we get 
\begin{equation}
 \boldsymbol\Upsilon(\widetilde{{\boldsymbol{\theta}}^g })^T IF_{i_0}(t_{i_0}, \widetilde{\boldsymbol{U}}_\tau, \underline{\mathbf{G}}) 
 = \boldsymbol{0}. 
 \label{EQ:8IF_eq2}
 \end{equation}
Solving (\ref{EQ:8IF_eq1}) and (\ref{EQ:8IF_eq2}), 
(as done for the i.i.d.~case in \cite{Ghosh:2015}),
we get a general expression for the IF of the RMDPDE as
\begin{eqnarray}
&&  IF_{i_0}(t_{i_0}, \widetilde{\boldsymbol{U}}_\tau, \underline{\mathbf{G}})
= \frac{1}{n}\boldsymbol{Q}(\widetilde{{\boldsymbol{\theta}}^g })^{-1} 
{\boldsymbol\Psi_n^{(0)}}(\widetilde{{\boldsymbol{\theta}}^g })^T 
\boldsymbol{D}_{\tau, i_0}^{(0)}(t_{i_0}; \widetilde{{\boldsymbol{\theta}}^g }),
\label{EQ:8IF_RMDPDE_gen_model}
\end{eqnarray} 
where $\boldsymbol{Q}({\boldsymbol{\theta}})=\left[{\boldsymbol\Psi_n^{(0)}}({\boldsymbol{\theta}})^T 
{\boldsymbol\Psi_n^{(0)}}({{\boldsymbol{\theta}}}) + 
\boldsymbol\Upsilon({{\boldsymbol{\theta}}}) \boldsymbol\Upsilon({{\boldsymbol{\theta}}})^T\right]$.
Clearly, this IF  is bounded in $t_{i_0}$ whenever 
$f_{i_0}(t_{i_0};\widetilde{{\boldsymbol{\theta}}}^g)^\tau \boldsymbol{u}_{i_0}^{(0)}(t_{i_0};\widetilde{{\boldsymbol{\theta}}}^g)$ 
is bounded, and this is the case at $\tau>0$ for most parametric models and common parametric restrictions.

If the contamination is in all the directions at the points 
$\mathbf{t} = (t_1, \cdots, t_n)$, the IF of the RMDPDE is given by 
\begin{eqnarray}
 IF_{o}(\mathbf{t}; \widetilde{\boldsymbol{U}}_\tau, \mathbf{\underline{G}}) 
 &=& \boldsymbol{Q}(\widetilde{{\boldsymbol{\theta}}^g })^{-1}  
 {\boldsymbol\Psi_n^{(0)}}(\widetilde{{\boldsymbol{\theta}}^g })^T
     \left[\frac{1}{n}\sum_{i=1}^n \boldsymbol{D}_{\tau, i}^{(0)}(t_{i}; \widetilde{{\boldsymbol{\theta}}^g })\right].\nonumber 
\end{eqnarray}

\subsection{Asymptotic Properties of the Proposed Test}
\label{SEC:asymp_composite}

Assume that $\Theta_0$ is a proper subset of the parameter space $\Theta$ 
which can be defined in terms of $r$ restrictions 
   $\boldsymbol\upsilon({\boldsymbol{\theta}}) = 0$
  such that the $p \times r$ matrix 
   $\boldsymbol\Upsilon({\boldsymbol{\theta}}) = \frac{\partial \boldsymbol\upsilon({\boldsymbol{\theta}})}{\partial {\boldsymbol{\theta}}}$
exists and is a continuous function of ${\boldsymbol{\theta}}$ with rank $r$. 

\begin{theorem}
	Suppose the model density satisfies the Lehmann and Ghosh-Basu conditions, $H_0$ is true with 
${\boldsymbol{\theta}}_0\in\Theta_0$ being the true parameter value, 
and (C1), (C2) and (C4) hold at ${\boldsymbol{\theta}}={\boldsymbol{\theta}}_0$. 
With
$
\widetilde{\boldsymbol\Sigma}_\tau({\boldsymbol{\theta}}_0) = [{\boldsymbol{J}_{\tau}}^{-1}({\boldsymbol{\theta}}_0)-{\boldsymbol{P}_{\tau}}({\boldsymbol{\theta}}_0)] {\boldsymbol{V}_{\tau}}({\boldsymbol{\theta}}_0)
[{\boldsymbol{J}_{\tau}}^{-1}({\boldsymbol{\theta}}_0)-{\boldsymbol{P}_{\tau}}({\boldsymbol{\theta}}_0)],
$
the asymptotic null distribution of the DPDTS$_C$ 	
$S_{\gamma}( {{\boldsymbol{\theta}}_n^\tau}, {\widetilde{{\boldsymbol{\theta}}}_n^\tau})$ is the distribution of 
$\sum_{i=1}^r ~  \widetilde{\zeta_i^{\gamma, \tau}}({\boldsymbol{\theta}}_0)Z_i^2,$
where 
$r = rank( {\boldsymbol{V}_{\tau}}({\boldsymbol{\theta}}_0) [ {\boldsymbol{J}_{\tau}}^{-1}({\boldsymbol{\theta}}_0) - {\boldsymbol{P}_{\tau}}({\boldsymbol{\theta}}_0)]{\boldsymbol{A}_{\gamma}}({\boldsymbol{\theta}}_0)
[{\boldsymbol{J}_{\tau}}^{-1}({\boldsymbol{\theta}}_0) - {\boldsymbol{P}_{\tau}}({\boldsymbol{\theta}}_0)] {\boldsymbol{V}_{\tau}}({\boldsymbol{\theta}}_0))$,
$Z_1, \cdots,Z_r$ are independent standard normals and
$\widetilde{\zeta_1^{\gamma, \tau}}({\boldsymbol{\theta}}_0)$, $\ldots$, $\widetilde{\zeta_r^{\gamma, \tau}}({\boldsymbol{\theta}}_0)$ 
are the nonzero eigenvalues of 
${\boldsymbol{A}_{\gamma}}({\boldsymbol{\theta}}_0)\widetilde{\boldsymbol\Sigma}_\tau({\boldsymbol{\theta}}_0)$.
\label{THM:8asymp_null_composite}
\end{theorem}

We can find approximate critical values of the asymptotic null distribution 
from the discussions in \cite{Basu/etc:2013a,Basu/etc:2013b}.
Next, we derive an asymptotic power approximation of the proposed DPDTS$_C$ 
at any point ${\boldsymbol{\theta}}^* \notin \Theta_0$, which can be used to determine minimum sample size 
requirement to attain any desired power as in the case of a simple hypothesis.
If ${\boldsymbol{\theta}}^* \notin \Theta_0$ is the true parameter value, then 
${\boldsymbol{\theta}}_n^\tau \displaystyle\mathop{\rightarrow}^\mathcal{P} {\boldsymbol{\theta}}^*$ and 
$\widetilde{{\boldsymbol{\theta}}}_n^\tau \displaystyle\mathop{\rightarrow}^\mathcal{P} {\boldsymbol{\theta}}_0$ 
for some ${\boldsymbol{\theta}}_0 \in \Theta_0$ and ${\boldsymbol{\theta}}^* \neq {\boldsymbol{\theta}}_0$. 
Then, assuming the Lehman conditions and Ghosh-Basu conditions along with (C1) and (C4) at 
${\boldsymbol{\theta}} = {\boldsymbol{\theta}}_0, {\boldsymbol{\theta}}^*$, we can show that 
  \begin{eqnarray}
  \sqrt{n}  \begin{pmatrix}
      ~~{\boldsymbol{\theta}}_n^\tau - {\boldsymbol{\theta}}^* ~\\
      ~~\widetilde{{\boldsymbol{\theta}}}_n^\tau - {\boldsymbol{\theta}}_0~
     \end{pmatrix} \rightarrow N_{2p}\left( \begin{bmatrix}
              ~\boldsymbol{0}~ \\
              ~\boldsymbol{0}~
             \end{bmatrix},  \begin{bmatrix}
                ~{\boldsymbol{\Sigma}_{\tau}}({\boldsymbol{\theta}}^*) & \boldsymbol{A}_{12}~\\ 
                ~\boldsymbol{A}_{12}^T & {\boldsymbol{P}_{\tau}}({\boldsymbol{\theta}}_0){\boldsymbol{V}_{\tau}}({\boldsymbol{\theta}}_0){\boldsymbol{P}_{\tau}}({\boldsymbol{\theta}}_0)~
                \end{bmatrix}\right),~~~\nonumber
                \label{EQ:8asymp_null_nonnull}
  \end{eqnarray}
for a $p\times p$ matrix $\boldsymbol{A}_{12}=\boldsymbol{A}_{12}({\boldsymbol{\theta}}^*,{\boldsymbol{\theta}}_0)$. 
Take 
$
\boldsymbol{M}_{1,\gamma}^{(i)}({\boldsymbol{\theta}}^*,{\boldsymbol{\theta}}_0) = 
\nabla d_\gamma(f_i(.;{\boldsymbol{\theta}}),f_i(.;{\boldsymbol{\theta}}_0))\big|_{{\boldsymbol{\theta}}={\boldsymbol{\theta}}^*}
$ and
$\boldsymbol{M}_{2,\gamma}^{(i)}({\boldsymbol{\theta}}^*,{\boldsymbol{\theta}}_0) = 
\nabla d_\gamma(f_i(.;{\boldsymbol{\theta}}^*),f_i(.;{\boldsymbol{\theta}}))\big|_{{\boldsymbol{\theta}}={\boldsymbol{\theta}}_0}
$.
\begin{itemize}
\item[(C5)]$\boldsymbol{M}_n^{j,\gamma}({\boldsymbol{\theta}}^*,{\boldsymbol{\theta}}_0) 
= n^{-1} \sum_{i=1}^n \boldsymbol{M}_{j,\gamma}^{(i)}({\boldsymbol{\theta}}^*,{\boldsymbol{\theta}}_0)
\rightarrow \boldsymbol{M}_{j,\gamma}({\boldsymbol{\theta}}^*,{\boldsymbol{\theta}}_0)$ element-wise as 
$n \rightarrow \infty$ for some $p$-vectors  $\boldsymbol{M}_{j,\gamma}$ ($j=1,2$).
\end{itemize}

\begin{theorem}\label{THM:8asymp_power_composite}
	Suppose the model density satisfies the Lehmann and Ghosh-Basu conditions and 
 ${\boldsymbol{\theta}}^*\notin \Theta_0$ for which (C1), (C4), and (C5) hold. 
Then, an approximation to the power function of the DPDTS$_C$ 
for testing (\ref{EQ:8composite_hypo}) at the significance level $\alpha$ is given by
\begin{eqnarray}
    \pi_{n,\alpha}^{\tau,\gamma} ({\boldsymbol{\theta}}^*) = 
  1 - \Phi \left( \frac{1}{\sqrt{n} \sigma_{\tau,\gamma}({\boldsymbol{\theta}}^*, {\boldsymbol{\theta}}_0)} 
\left(\frac{s_\alpha^{\tau,\gamma}}{2}-\sum_{i=1}^n d_\gamma(f_i(.;{\boldsymbol{\theta}}^*),f_i(.;{\boldsymbol{\theta}}_0))\right)\right),\nonumber
\end{eqnarray}
where $s_\alpha^{\tau,\gamma}$ is $(1-\alpha)$-th quantile of the asymptotic null distribution of
$S_{\gamma}( {{\boldsymbol{\theta}}_n^\tau}, \widetilde{{\boldsymbol{\theta}}}_n^\tau)$, 
\begin{eqnarray}
    \sigma_{\tau,\gamma}^2({\boldsymbol{\theta}}^*,{\boldsymbol{\theta}}_0) 
    = \boldsymbol{M}_{1,\gamma}^T {\boldsymbol{\Sigma}_{\tau}} \boldsymbol{M}_{1,\gamma} 
    + \boldsymbol{M}_{1,\gamma}^T \boldsymbol{A}_{12} \boldsymbol{M}_{2,\gamma} 
+ \boldsymbol{M}_{2,\gamma}^T \boldsymbol{A}_{12}^T \boldsymbol{M}_{1,\gamma} 
    + \boldsymbol{M}_{2,\gamma}^T {\boldsymbol{P}_{\tau}} {\boldsymbol{V}_{\tau}} {\boldsymbol{P}_{\tau}} \boldsymbol{M}_{2,\gamma}.\nonumber
\end{eqnarray}
\end{theorem}

\begin{corollary}
For ${\boldsymbol{\theta}}^* \ne {\boldsymbol{\theta}}_0$, the probability of rejecting $H_0$ 
in (\ref{EQ:8composite_hypo}) at level $\alpha > 0 $ based on the DPDTS$_C$ 
tends to $1$ as $n \rightarrow \infty$, provided 
$\frac{1}{n} \displaystyle\sum_{i=1}^n d_\gamma(f_i(.;{\boldsymbol{\theta}}^*),f_i(.;{\boldsymbol{\theta}}_0)) = O(1). $ 
\end{corollary}

\subsection{Robustness Properties of the Test}
\label{SEC:8robust_composite}


Using the functional form of ${{{\boldsymbol{\theta}}}_n^\tau}$ and ${\widetilde{{\boldsymbol{\theta}}}_n^\tau}$  
and ignoring the multiplier $2$ in our test statistic, 
the functional corresponding to the DPDTS$_C$ is
\begin{eqnarray}\label{EQ:8DPDTS2_IF}
    S_{\gamma, \tau}^{(1)}(\underline{\mathbf{G}}) = \sum_{i=1}^n 
    d_\gamma(f_i(.;\boldsymbol{U}_\tau(\underline{\mathbf{G}})),f_i(.;\widetilde{\boldsymbol{U}}_\tau(\underline{\mathbf{G}}))).
    \nonumber
\end{eqnarray}
Clearly, this depends on the sample size $n$, implying the same dependency in its IF. 
Consider the contaminated distribution  $G_{i,\epsilon}$ defined in Section
\ref{SEC:8IF_simple_test} and assume contamination in only one fixed direction-$i_0$. 
The first order IF of $S_{\gamma, \tau}^{(1)}(\underline{\mathbf{G}})$ is
\vspace{-0.1in}
\begin{eqnarray} 
IF_{i_0}(t_{i_0}, S_{\gamma, \tau}^{(1)}, \underline{\mathbf{G}})  = \frac{\partial}{\partial\epsilon}
 S_{\gamma, \tau}^{(1)}(\underline{\mathbf{G}}_{i_0,\epsilon}) \big|_{\epsilon=0} 
&=& 
n \boldsymbol{M}_n^{1,\gamma}(\boldsymbol{U}_\tau(\underline{\mathbf{G}}), \widetilde{\boldsymbol{U}}_\tau(\underline{\mathbf{G}}))^T 
 IF_{i_0}(t_{i_0}, \boldsymbol{U}_\tau, \underline{\mathbf{G}}) \nonumber 
 \\&+& 
n \boldsymbol{M}_n^{2,\gamma}({\boldsymbol{U}}_\tau(\underline{\mathbf{G}}), \widetilde{\boldsymbol{U}}_\tau(\underline{\mathbf{G}}))^T 
 IF_{i_0}(t_{i_0}, \widetilde{\boldsymbol{U}}_\tau, \underline{\mathbf{G}}),\nonumber
\end{eqnarray}
where $IF_{i_0}(t_{i_0}, \widetilde{\boldsymbol{U}}_\tau, \underline{\mathbf{G}})$ is the IF
of the RMDPD functional $\widetilde{\boldsymbol{U}}_{\tau}$ under $H_0$. 
If the null hypothesis is true with 
$\underline{\mathbf{G}}=\underline{\mathbf{F}}_{{\boldsymbol{\theta}}_0}$ 
for some ${\boldsymbol{\theta}}_0 \in \Theta_0$, then 
$\boldsymbol{U}_\tau(\underline{\mathbf{F}}_{{\boldsymbol{\theta}}_0}) 
=\widetilde{\boldsymbol{U}}_\tau(\underline{\mathbf{F}}_{{\boldsymbol{\theta}}_0}) = {\boldsymbol{\theta}}_0$ and 
$\boldsymbol{M}_{j,\gamma}^{(i)}({\boldsymbol{\theta}}_0,{\boldsymbol{\theta}}_0)=\boldsymbol{0}$ for $j=1, 2$. 
Hence Hampel's first-order IF of the DPDTS$_C$  
is again zero at the composite null.

Similarly, at $\underline{\mathbf{G}}=\underline{\mathbf{F}}_{{\boldsymbol{\theta}}_0}$, 
the second order IF of the DPDTS$_C$ functional $ S_{\gamma, \tau}^{(1)}$ is
%
\begin{align}
IF_{i_0}^{(2)}(t_{i_0}, S_{\gamma, \tau}^{(1)}, \underline{\mathbf{F}}_{{\boldsymbol{\theta}}_0})
= n  \boldsymbol{D}_{\tau,i_0}(t_{i_0}, {\boldsymbol{\theta}}_0)^T \boldsymbol{A}_n^\gamma({\boldsymbol{\theta}}_0) 
\boldsymbol{D}_{\tau,i_0}(t_{i_0}, {\boldsymbol{\theta}}_0),  
\nonumber 
\end{align}
where 
$
\boldsymbol{D}_{\tau,i_0}(t_{i_0}, {\boldsymbol{\theta}}_0)
=\left[IF_{i_0}(t_{i_0}, \boldsymbol{U}_\tau, \underline{\mathbf{F}}_{{\boldsymbol{\theta}}_0}) - 
IF_{i_0}(t_{i_0}, \widetilde{\boldsymbol{U}}_\tau, \underline{\mathbf{F}}_{{\boldsymbol{\theta}}_0})\right].
$
Clearly, this IF is bounded 
if  the corresponding MDPDEs over $\Theta_0$ and $\Theta$ both have bounded IFs. 
However, the boundedness of the IF of 
the MDPDE over $\Theta$ implies the same under any restricted subspace $\Theta_0$
and this holds for most parametric models if $\tau>0$, but the IF is unbounded at $\tau=0$.

Next, considering the contamination in all directions at $\mathbf{t} = (t_1, \ldots, t_n)$, 
the first order IF of the proposed DPDTS$_C$ is again zero at any point inside $\Theta_0$
and its second order IF at the null is given by
\begin{eqnarray}
IF_o^{(2)}(\mathbf{t}, T_{\gamma, \tau}^{(1)}, \underline{\mathbf{F}}_{{\boldsymbol{\theta}}_0}) 
= n \boldsymbol{D}_{\tau,o}(\mathbf{t}, {\boldsymbol{\theta}}_0)^T \boldsymbol{A}_n^\gamma({\boldsymbol{\theta}}_0) 
\boldsymbol{D}_{\tau,o}(\mathbf{t}, {\boldsymbol{\theta}}_0),\nonumber
\end{eqnarray}
where $\boldsymbol{D}_{\tau,o}(\mathbf{t}, {\boldsymbol{\theta}}_0)
= \left[IF_o(\mathbf{t}, \boldsymbol{U}_\tau, \underline{\mathbf{F}}_{{\boldsymbol{\theta}}_0}) - 
IF_o(\mathbf{t}, \widetilde{\boldsymbol{U}}_\tau, \underline{\mathbf{F}}_{{\boldsymbol{\theta}}_0})\right]$.
This implies the robustness for $\tau >0$.

%


The level and power influence functions  of the proposed test for this case are similar to 
that in the simple hypothesis case (Section \ref{SEC:8IF_power_simple_test}).

\section{Application (I): Normal Linear Regression}\label{SEC:linear_regression}

Consider the linear regression model
$y_i = {\boldsymbol{x}_i}^T{{\mathbf {\boldsymbol{\beta}}}} + \epsilon_i$,  for $i = 1, \ldots, n$,
%
where the error $\epsilon_i$'s are assumed to be i.i.d.~normal with mean
zero and variance $\sigma^2$; ${{\mathbf {\boldsymbol{\beta}}}}= ({{\beta}}_1, \ldots, {{\beta}}_p)^T$ and
${\boldsymbol{x}_i}^T = (x_{i1}, \ldots, x_{i,p})$   denote the regression coefficients  and
the $i$-th observation for the covariates, respectively. 
Here, we  assume ${\boldsymbol{x}_i}$ to be fixed so that 
$y_i \sim N({\boldsymbol{x}_i}^T{{\mathbf {\boldsymbol{\beta}}}}, \sigma^2)$ for each $i$. 
Clearly the $y_i$'s are independent but not identically distributed.  
The MDPDEs of $\boldsymbol{\beta}$ and $\sigma^2$ and their properties are described in 
Section S2.1 of the Online Supplement.

\subsection{Testing for the Regression Coefficients with Known $\sigma$}

The simple hypothesis on the regression coefficient ${\boldsymbol{\beta}}(={\boldsymbol{\theta}})$ 
assuming the error variance $\sigma^2$ to be known, say $\sigma^2 =\sigma_0^2$, is 
\begin{equation}
H_0 : {\boldsymbol{\beta}}={\boldsymbol{\beta}}_0,  ~~~~\mbox{ against }~~~~  H_1 : {\boldsymbol{\beta}}\neq{\boldsymbol{\beta}}_0,
\label{EQ:beta_simple_Hyp}
\end{equation}
for some pre-specified ${\boldsymbol{\beta}}_0(={\boldsymbol{\theta}}_0)$.  

Consider the test statistics
$T_{\gamma}( {{\boldsymbol{\beta}}_n^\tau}, {{\boldsymbol{\beta}}_0})$  for testing (\ref{EQ:beta_simple_Hyp}),
where ${\boldsymbol{\beta}}_n^\tau$ is the MDPDE of ${\boldsymbol{\beta}}$ with tuning parameter $\tau$  and known $\sigma=\sigma_0$. 
Using the normal density, we get
\begin{eqnarray}
T_{\gamma}( {{\boldsymbol{\beta}}_n^\tau}, {{\boldsymbol{\beta}}_0})  &=& 
\frac{2\sqrt{1+\gamma}}{\gamma(\sqrt{2\pi}\sigma_0)^\gamma} \left[n-\sum_{i=1}^n~
e^{-\frac{\gamma({\boldsymbol{\beta}}_n^\tau-{\boldsymbol{\beta}}_0)^T({\boldsymbol{x}_i}{\boldsymbol{x}_i}^T)({\boldsymbol{\beta}}_n^\tau-{\boldsymbol{\beta}}_0)}{2(\gamma(\sigma_n^\tau)^2 
+ \sigma_0^2)}} \right],  \mbox{ if } \gamma>0,\nonumber 
\\
\mbox{and }  T_{0}( {{\boldsymbol{\beta}}_n^\tau}, {{\boldsymbol{\beta}}_0})  &=&  
\frac{({\boldsymbol{\beta}}_n^\tau-{\boldsymbol{\beta}}_0)^T (\boldsymbol{X}^T\boldsymbol{X}) ({\boldsymbol{\beta}}_n^\tau-{\boldsymbol{\beta}}_0)}{\sigma_0^2},\nonumber
\end{eqnarray}
where $\boldsymbol{X}=[\boldsymbol{x}_1 \cdots \boldsymbol{x}_n]^T$.
The estimator ${\boldsymbol{\beta}}_n^{(0)}$, the MDPDE with $\tau=0$, is indeed the MLE of ${\boldsymbol{\beta}}$. 
The usual LRT statistics for this problem is $-2\log 
\left[\frac{\prod_{i=1}^n \phi(y_i, {\boldsymbol{x}_i}^T{\boldsymbol{\beta}}_0, \sigma_0)}{
\prod_{i=1}^n \phi(y_i; {\boldsymbol{x}_i}^T{\boldsymbol{\beta}}_n^{(0)}, \sigma_0)}\right]$;
after simplification, this statistic is
the same as $T_{0}( {{\boldsymbol{\beta}}_n^{(0)}}, {{\boldsymbol{\beta}}_0})$. 
Hence the proposed test is nothing but a robust generalization of the likelihood ratio test.
Here $\phi(\cdot, \mu, \sigma)$ refers to the $N(\mu, \sigma^2)$ density.

\subsubsection{Asymptotic Properties}

Suppose that the conditions (R1)--(R2) of \cite{Ghosh/Basu:2013}
hold true and assume 
\begin{itemize}
\item[(C6)] $\frac{1}{n} (\boldsymbol{X}^T\boldsymbol{X})$ converges point-wise to some positive definite matrix 
$\boldsymbol\Sigma_x$  as $n\rightarrow\infty$.
\end{itemize}
 Then, the corresponding limiting matrices simplify to 
 $ {\boldsymbol{J}_{\tau}}({\boldsymbol{\beta}}_0) =  \zeta_\tau \boldsymbol\Sigma_x$, 
 ${\boldsymbol{V}_{\tau}}({\boldsymbol{\beta}}_0) =  \zeta_{2\tau} \boldsymbol\Sigma_x$ and 
 ${\boldsymbol{A}_{\gamma}}({\boldsymbol{\beta}}_0) =  (1+\gamma)\zeta_\gamma \boldsymbol\Sigma_x$, 
where $\zeta_\tau= (2\pi)^{-\frac{\tau}{2}} \sigma^{-(\tau+2)}(1+\tau)^{-\frac{3}{2}}$.

Theorem \ref{THM:8asymp_null_simple_test} gives the asymptotic null distribution of 
$T_{\gamma}( {{\boldsymbol{\beta}}_n^\tau}, {{\boldsymbol{\beta}}_0}) $ under $H_0 : {\boldsymbol{\beta}}={\boldsymbol{\beta}}_0$, 
which turns out to be a scalar multiple of a $\chi_p^2$ distribution 
with the multiplier 
$\zeta_1^{\gamma, \tau} =  
(\sqrt{2\pi}\sigma_0)^{-\gamma} (1+\gamma)^{-\frac{1}{2}} 
\left(1 + \frac{\tau^2}{1+2\tau}\right)^{\frac{3}{2}}. $ 
The critical region for testing (\ref{EQ:beta_simple_Hyp}) at the significance level $\alpha$ is given by 
$
\left\{T_{\gamma}( {{\boldsymbol{\beta}}_n^\tau}, {{\boldsymbol{\beta}}_0}) > \zeta_1^{\gamma, \tau}\chi^2_{p, \alpha}\right\},
$
where $\chi^2_{p, \alpha}$ is the $(1-\alpha)$-th quantile of the $\chi_p^2$ distribution. 
At $\gamma=\tau=0$, we have $\zeta_1^{0, 0}=1$ so that $T_{0}( {{\boldsymbol{\theta}}_n^{(0)}}, {{\boldsymbol{\theta}}_0})$ 
asymptotically follows a $\chi_p^2$ distribution under $H_0$, 
 as expected from its relation to the LRT.

\begin{figure}[!h]
	\centering
	\subfloat[$p=2$]{
		\includegraphics[width=0.33\textwidth]{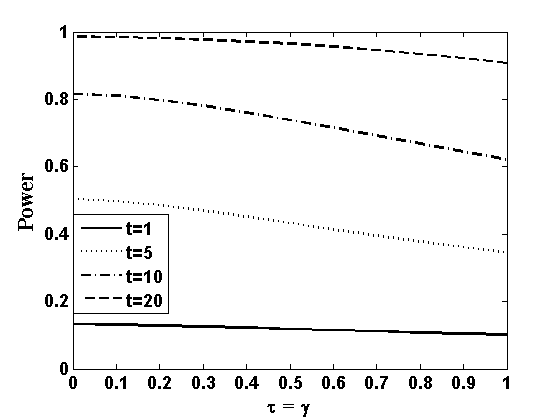}
		\label{FIG:9simple_contPower_p2}}
	~ 
	\subfloat[$p=10$]{
		\includegraphics[width=0.33\textwidth]{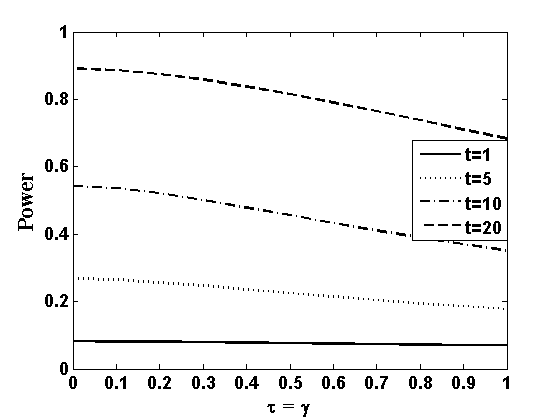}
		\label{FIG:9simple_contPower_p10}}
	\\
	\subfloat[$p=50$]{
		\includegraphics[width=0.33\textwidth]{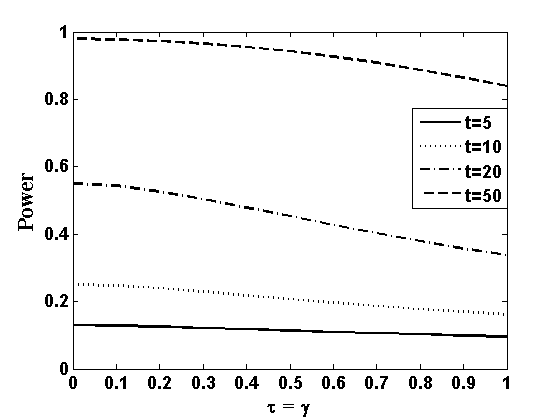}
		\label{FIG:9simple_contPower_p50}}
	~ 
	\subfloat[$p=200$]{
		\includegraphics[width=0.33\textwidth]{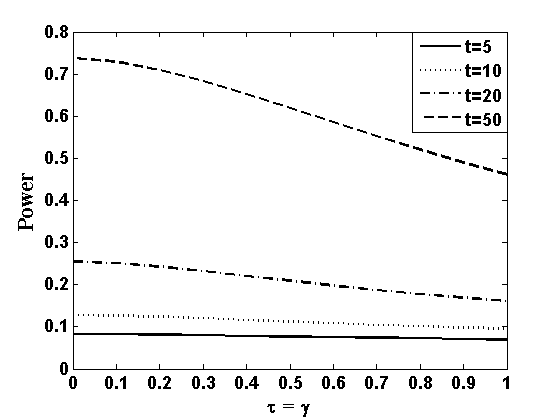}
		\label{FIG:9simple_contPower_p200}}
	\caption{Asymptotic contiguous power of the simple DPD based test of ${\boldsymbol{\beta}}$ for different values of 
		$ t = \boldsymbol\Delta^T\boldsymbol\Sigma_x\boldsymbol\Delta$ and $p$, the number of explanatory variables}
	\label{FIG:9simple_contPower}
\end{figure}

Next we study the asymptotic power of the proposed test.
We derive its asymptotic power under the contiguous alternatives 
$H_{1,n}$ using Corollary \ref{COR:8contiguous_power_simple_test}. 
The asymptotic distribution of  $T_{\gamma}( {{\boldsymbol{\beta}}_n^\tau}, {{\boldsymbol{\beta}}_0}) $
under  $H_{1,n}$ is $\zeta_1^{\gamma, \tau}\chi^2_{p,\delta}$ with
$\delta = \frac{1}{\upsilon_\tau^{\boldsymbol{\beta}}}  \boldsymbol\Delta^T\boldsymbol\Sigma_x\boldsymbol\Delta$,
where $\upsilon_\tau^{{\mathbf {\boldsymbol{\beta}}}}= \sigma^2 \left(1 + \frac{\tau^2}{1+2\tau}\right)^{\frac{3}{2}}$.
Thus its asymptotic contiguous power is
$
P_{\tau, \gamma}(\boldsymbol\Delta,0; \alpha) $ 
$= P\left(\zeta_1^{\gamma, \tau}W_{p,\delta} > \zeta_1^{\gamma, \tau}\chi^2_{p, \alpha}\right)$
$=1 - G_{p,\delta}(\chi^2_{p, \alpha}),
$
where  $G_{p,\delta}$ is the distribution function of $\chi^2_{p,\delta}$. 
Figure \ref{FIG:9simple_contPower} shows the asymptotic power over the tuning parameters 
$\gamma=\tau$ for different values of $\boldsymbol\Delta^T\boldsymbol\Sigma_x\boldsymbol\Delta~(= t,~\mbox{ say})$; 
it does not depend on the tuning parameter $\gamma$. 
The power depends on the distance ($\boldsymbol\Delta$) of the 
contiguous alternatives from null and the limiting second order moments ($\boldsymbol\Sigma_x$) of the covariates through 
$ t = \boldsymbol\Delta^T\boldsymbol\Sigma_x\boldsymbol\Delta$; for any fixed $\tau=\gamma$ it increases as the 
value of $t$ increases. It also depends on the number ($p$) of explanatory variables used.
In Figure \ref{FIG:9simple_contPower}, we show the cases of $p$ is 2 and 10 
as well as $p=50$ and 200. 
For fixed values of $\gamma=\tau$ and $\boldsymbol\Delta_1^T\boldsymbol\Sigma_x\boldsymbol\Delta_1$, 
the power decreases as $p$ increases; this is expected as 
the number of components of ${\boldsymbol{\beta}}$ increases with $p$.
The power against any contiguous alternative and any model 
is seen to decrease slightly with increasing  $\tau$
which brings in 
the non-centrality parameter $\delta$
and hence the asymptotic variance  $\upsilon_\tau^{\boldsymbol{\beta}}$  of each element of 
$(\boldsymbol{X}^T\boldsymbol{X})^{1/2}\boldsymbol{\beta}_n^\tau$. 
As $\upsilon_\tau^{\boldsymbol{\beta}}$ increases slightly with $\tau$,
the efficiency of the MDPDE and the asymptotic contiguous power of the DPDTS decrease slightly.

\begin{figure}[!th]
	\centering
	\subfloat[$\boldsymbol\Delta^T\boldsymbol\Sigma_x\boldsymbol\Delta=5$]{
		\includegraphics[width=0.33\textwidth]{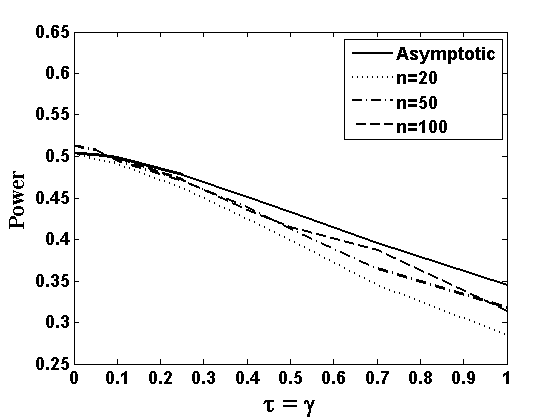}
		\label{FIG:9simple_contPower_p2_comp}}
	~ 
	\subfloat[$\boldsymbol\Delta^T\boldsymbol\Sigma_x\boldsymbol\Delta=20$]{
		\includegraphics[width=0.33\textwidth]{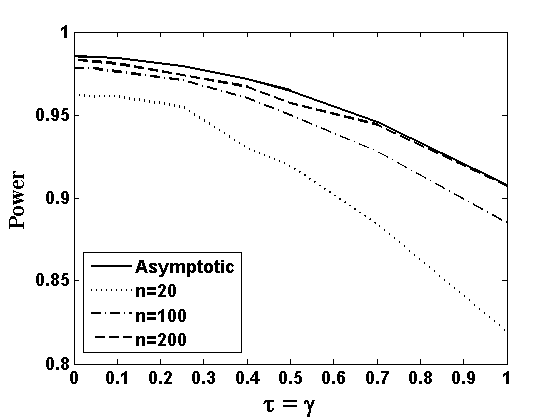}
		\label{FIG:9simple_contPower_p10_comp}}
	\caption{Comparison of finite-sample empirical power at different sample sizes $n$ with 
		asymptotic contiguous power for the simple DPDTS of ${\boldsymbol{\beta}}$ for $p=2$ and $\sigma_x^2=5$}
	\label{FIG:9simple_contPower_comp}
\end{figure}

We simulated finite-sample situations
with different sample sizes $n$ and values of $ t = \boldsymbol\Delta^T\boldsymbol\Sigma_x\boldsymbol\Delta$ and $p$. 
The convergence of the finite-sample power to the asymptotic value depends on 
the convergence rate of $\frac{1}{n} (\boldsymbol{X}^T\boldsymbol{X})$ in Condition (C6);
we chose  $X$ to be a $p$-variate normal distribution with mean $0$ and covariance matrix $\sigma_x^2\boldsymbol{I}_p$
so that (C6) holds with $\boldsymbol\Sigma_x= \sigma_x^2\boldsymbol{I}_p$.
Results from two such simulations with $p=2$ and $\sigma_x^2=5$
are presented in Figure \ref{FIG:9simple_contPower_comp}; other cases have similar patterns. 
Clearly the finite-sample powers are close to the asymptotic power in moderate sample sizes like $n=100$
and the convergence rate is little slower for larger $ t = \boldsymbol\Delta^T\boldsymbol\Sigma_x\boldsymbol\Delta$.

\subsubsection{Robustness Results}\label{SEC:9robust_test_knownSigma}

As the first order IF of DPDTS $T_{\gamma}( {{\boldsymbol{\beta}}_n^\tau}, {{\boldsymbol{\beta}}_0})$ 
is zero at any simple null, we measure the stability of the proposed test by the second order IF. 
For contamination in a single direction ($i_0^{\rm th}$ direction), 
the second order IF at the null  ${\boldsymbol{\beta}}={\boldsymbol{\beta}}_0$ simplifies to
\begin{eqnarray}
&& IF_{i_0}^{(2)}(t_{i_0}, T_{\gamma, \tau}^{(1)}, \mathbf{\underline{F}}_{{\boldsymbol{\theta}}_0}) 
\nonumber \\&& ~~
= (1+\gamma)\zeta_\gamma(1+\tau)^{3}n[{\boldsymbol{x}_{i_0}}^T(\boldsymbol{X}^T\boldsymbol{X})^{-1}{\boldsymbol{x}_{i_0}}] (t_{i_0} - {\boldsymbol{x}_{i_0}}^T{{\mathbf {\boldsymbol{\beta}}_0}})^2
e^{-\frac{\tau(t_{i_0} - {\boldsymbol{x}_{i_0}}^T{{\mathbf {\boldsymbol{\beta}}_0}})^2}{\sigma_0^2}}.\nonumber
\end{eqnarray}   
As expected, the IF depends on the outliers and the leverage points through 
$(t_{i_0} - {\boldsymbol{x}_{i_0}}^T{{\mathbf {\boldsymbol{\beta}}_0}})$ and $[{\boldsymbol{x}_{i_0}}^T(\boldsymbol{X}^T\boldsymbol{X})^{-1}{\boldsymbol{x}_{i_0}}]$.
The LIF and PIF of the proposed DPDTS under contiguous alternatives 
is presented in Section S2.2 of the Online Supplement. 

Both the IF and the PIF are bounded with respect to the contamination point $t_{i_0}$ for any $\tau>0$ 
implying their stability against contamination. 
Both are unbounded at $t_{i_0}$ for the LRT at $\gamma=\tau=0$
indicating its non-robustness. 
The LIF of this test is identically zero for all $\tau, \gamma \geq 0$
implying no asymptotic influence of contiguous contamination on its size.

\section{Application (II): Generalized Linear Model}\label{SEC:glm}

In a generalized linear model (GLM), the response variables $Y_i$ are independent 
and have an exponential family distribution with density
$f(y_i;{{\theta}}_i,\phi) = \exp\left\{ \frac{y_i{{\theta}}_i - b({{\theta}}_i)}{a(\phi)} + c(y_i,\phi)  \right\}$;
the canonical parameter ${{\theta}}_i$ depends on the predictor ${\boldsymbol{x}_i}$ and 
$\phi$ is a nuisance scale parameter.
The mean $\mu_i$ of $Y_i$ satisfies $g(\mu_i) = \eta_i = {\boldsymbol{x}_i}^T{\boldsymbol{\beta}},$ for a monotone differentiable link function $g$ 
and linear predictor $\eta_i = {\boldsymbol{x}_i}^T{\boldsymbol{\beta}}$.

The  GLMs with fixed predictors fit the general  I-NH set-up.
The properties of the MDPDEs of ${\boldsymbol{\theta}} = ({\boldsymbol{\beta}},\phi)$ 
in the GLM were derived in \cite{Ghosh/Basu:2014}
and are presented in Section S3 of the Online Supplement.

Suppose we have a sample of size $n$ from a GLM with parameter 
${\boldsymbol{\theta}}=({\boldsymbol{\beta}},~\phi) \in \Theta=\mathbb{R}^p\times [0, \infty)$  
and want to test the hypothesis 
\begin{eqnarray}
H_0 : \boldsymbol{L}^T{\boldsymbol{\beta}} = \boldsymbol{l}_0 ~~~~\mbox{ against }~~~~~ 
H_1 : \boldsymbol{L}^T{\boldsymbol{\beta}} \neq \boldsymbol{l}_0, ~~
\label{EQ:9Gen_lin_hypothesis}
\end{eqnarray}  
where $\boldsymbol{L}$ is a $p\times r$ known matrix with $p\geq r$ and 
$\boldsymbol{l}_0$ is a real $r$-vector.
We assume that the nuisance parameter  $\phi$ is unknown.

The DPD based test statistic (DPDTS$_C$) for testing this problem is  
$$
S_{\gamma}( {{\boldsymbol{\theta}}_n^\tau}, {\widetilde{{\boldsymbol{\theta}}}_n^\tau}) 
= 2 \sum_{i=1}^n d_\gamma(f_i(.;(\widehat{\boldsymbol{\beta}}_n^\tau, \widehat\phi_n^\tau)),
f_i(.;(\widetilde{{\boldsymbol{\beta}}_n}^\tau, \widetilde{\phi_n}^\tau))),
$$
where ${{\boldsymbol{\theta}}_n^\tau}=(\widehat{\boldsymbol{\beta}}_n^\tau, \widehat\phi_n^\tau)$ is the unrestricted MDPDE,
${\widetilde{{\boldsymbol{\theta}}}_n^\tau}=(\widetilde{{\boldsymbol{\beta}}_n}^\tau, \widetilde{\phi_n}^\tau)$ is the restricted MDPDE 
under $H_0$ corresponding to the tuning parameter $\tau$.

To derive the asymptotic distribution of the RMDPDE $(\widetilde{{\boldsymbol{\beta}}_n}^\tau, \widetilde{\phi_n}^\tau)$ 
of $({\boldsymbol{\beta}}, \phi)$ from Theorem \ref{THM:8asymp_RMDPDE}, some simple matrix algebra leads to 
\begin{eqnarray}
\boldsymbol{P}_n^\tau({\boldsymbol{\beta}}, \sigma) 
= n~\begin{bmatrix}
~~\boldsymbol\Psi_{n,11.2}^{-1}\left[\boldsymbol{I}_p - \boldsymbol{L}
\{\boldsymbol{L}^T\boldsymbol\Psi_{n,11.2}^{-1}\boldsymbol{L}\}^{-1}\boldsymbol{L}^T\boldsymbol\Psi_{n,11.2}^{-1}\right] 
& - \boldsymbol{M}_{11}\boldsymbol{X}^T\boldsymbol\Gamma_{12}^{(\tau)}\mathbf{1}\boldsymbol\Psi_{n,22.1}^{-1} ~~\\
~~-\boldsymbol\Psi_{n,22.1}^{-1}\mathbf{1}^T\boldsymbol\Gamma_{12}^{(\tau)}\boldsymbol{X} \boldsymbol{M}_{11}
& \boldsymbol\Psi_{n,22.1}^{-1} ~~
\end{bmatrix},\nonumber
\end{eqnarray}
where, for any $i,j=1,2$, 
$\boldsymbol\Psi_{n,ii.j} = \boldsymbol{X}^T\boldsymbol\Gamma_{jj}^{(\tau)}\boldsymbol{X} 
- \boldsymbol{X}^T\boldsymbol\Gamma_{ij}^{(\tau)}\mathbf{1} 
 (\mathbf{1}^T\boldsymbol\Gamma_{jj}^{(\tau)}\mathbf{1})^{-1} \mathbf{1}^T\boldsymbol\Gamma_{ji}^{(\tau)}\boldsymbol{X}$ 
with $\boldsymbol\Gamma_{ij}^{(\tau)}$  defined in Section 4 of the online Supplement,
and $\boldsymbol{M}_{11} = (\boldsymbol{X}^T\boldsymbol\Gamma_{11}^{(\tau)}\boldsymbol{X})^{-1} $.

\begin{corollary}\label{THM:10asymp_RMDPDE_GLM}
	Suppose the Ghosh-Basu conditions hold with respect to $\Theta_0$. 
	The RMDPDE $(\widetilde{{\boldsymbol{\beta}}_n}, \widetilde{\phi_n})$  exists 
	and is consistent for ${\boldsymbol{\theta}}_0 = ({\boldsymbol{\beta}}^g,\phi^g)$,
	true parameter value under $\Theta_0$. 
	The asymptotic distribution of 
	$\boldsymbol\Omega_n^{-\frac{1}{2}}	\boldsymbol{P}_n 
	\left[\sqrt n \left((\widetilde{{\boldsymbol{\beta}}_n}, \widetilde{\phi_n}) 
	- (\widetilde{{\boldsymbol{\beta}}^g}, \widetilde{\phi^g})\right)\right]$ is $(p+1)$-dimensional normal with 
	mean $\boldsymbol{0}$ and variance $\boldsymbol{I}_{p+1}$, 
	where $\boldsymbol{P}_n=\boldsymbol{P}_n^\tau(\widetilde{{\boldsymbol{\beta}}^g}, \widetilde{\phi^g})$ and 
	$\boldsymbol\Omega_n=\boldsymbol\Omega_n^\tau(\widetilde{{\boldsymbol{\beta}}^g}, \widetilde{\phi^g})$, 
	with $\boldsymbol\Omega_n({\boldsymbol{\beta}},\phi)$ defined in Section S3 of the Online Supplement.
\end{corollary}

As in the unrestricted case, 
the restricted MDPDE of ${\boldsymbol{\beta}}$ and $\phi$ may not be asymptotically independent. 
They are independent if $\gamma_{12i}^{1+2\tau}=0, \gamma_{1i}^{1+\tau}\gamma_{2i}^{1+\tau}=0, ~\forall i$.

To derive asymptotic distribution of the DPDTS$_C$, we assume 
fixed covariates ${\boldsymbol{x}_i}$'s for which the matrices
$\boldsymbol\Psi_n^\tau(\widetilde{{\boldsymbol{\theta}}}^g) $ and 
$\boldsymbol\Omega_n^\tau(\widetilde{{\boldsymbol{\theta}}}^g)$, defined in Section S3 of the Online Supplement, 
converge element-wise as $n \rightarrow \infty$ to some $(p+1) \times (p+1)$ invertible matrices ${\boldsymbol{J}_{\tau}}$ and ${\boldsymbol{V}_{\tau}}$, respectively.
Consider the partitions 
\begin{equation}
{\boldsymbol{J}_{\tau}}({\boldsymbol{\beta}}, \sigma) = \begin{bmatrix}
~~\boldsymbol{J}_{11,\tau} & \boldsymbol{J}_{12,\tau}~\\
~~\boldsymbol{J}_{12,\tau}^T & \boldsymbol{J}_{22,\tau}~
\end{bmatrix},   ~~~~~ \mbox{ and } ~~~~~
{\boldsymbol{V}_{\tau}}({\boldsymbol{\beta}}, \sigma) = \begin{bmatrix}
~~\boldsymbol{V}_{11,\tau} & \boldsymbol{V}_{12,\tau}~\\
~~\boldsymbol{V}_{12,\tau}^T & \boldsymbol{V}_{22,\tau}~
\end{bmatrix},\nonumber
\end{equation}
where $\boldsymbol{J}_{11,\tau}$ and $\boldsymbol{V}_{11,\tau}$ are of order $p\times p$.
We suppress $\tau$ in above notation whenever it is clear from the context.
Then, the asymptotic null distribution of the DPDTS$_C$ 
$S_{\gamma}( {{\boldsymbol{\theta}}_n^\tau}, {\widetilde{{\boldsymbol{\theta}}}_n^\tau})$ 
follows directly from Theorem \ref{THM:8asymp_null_composite} 
provided the Ghosh-Basu conditions hold for the model under $H_0$.

\begin{corollary}
For the GLM set-up, assume that its density satisfies the Lehmann 
and Ghosh-Basu conditions under $\Theta_0$. Then the asymptotic null distribution of the 
DPDTS$_C$ 	$S_{\gamma}( {{\boldsymbol{\theta}}_n^\tau}, {\widetilde{{\boldsymbol{\theta}}}_n^\tau})$ 
is the same as that of $\sum_{i=1}^r ~  \zeta_i^{\gamma, \tau}({\boldsymbol{\theta}}_0)Z_i^2,$
where $Z_1, \cdots,Z_r$ are independent standard normal variables, and
$\zeta_1^{\gamma, \tau}({\boldsymbol{\theta}}_0)$, $\cdots, \zeta_r^{\gamma, \tau}({\boldsymbol{\theta}}_0)$ are 
$r$ nonzero eigenvalues of the matrix
$$
\boldsymbol{E} = 
(1+\gamma)\boldsymbol{J}_{11,\gamma}\boldsymbol{J}_{11.2}^{-1}\boldsymbol{L}\boldsymbol{N}_{11}\boldsymbol{L}^T
\boldsymbol{J}_{11.2}^{-1}\boldsymbol{V}_{11}\boldsymbol{J}_{11.2}^{-1}\boldsymbol{L}
\boldsymbol{N}_{11}\boldsymbol{L}^T\boldsymbol{J}_{11.2}^{-1},
$$
where $\boldsymbol{J}_{ii.j} = \boldsymbol{J}_{ii,\tau} - \boldsymbol{J}_{ij,\tau}\boldsymbol{J}_{jj,\tau}^{-1}\boldsymbol{J}_{ji,\tau}^T$ 
for $i,j=1,2; i\neq j$ and $\boldsymbol{N}_{11}=(\boldsymbol{L}^T\boldsymbol{J}_{11.2}^{-1}\boldsymbol{L})^{-1}$. 
\label{THM:10asymp_null_GLM_test}
\end{corollary}

This result can be used to obtain the critical values of the proposed DPD based test.
The other asymptotic results regarding power and robustness of the test can be derived 
by direct application of the general theory of Section \ref{SEC:8composite_testing}.
For instance, the second order IF
of the test statistics at the null hypothesis, when there is contamination in only one fixed direction-$i_0$, is
\begin{eqnarray}	
IF_{i_0}^{(2)}(t_{i_0}, S_{\gamma, \tau}^{(1)}, \underline{\mathbf{F}}_{{\boldsymbol{\theta}}_0})
= n(1+\gamma) \cdot \boldsymbol{W}^T \boldsymbol\Psi_n^\gamma \boldsymbol{W},  
\label{EQ:10IF_DPD_test_GLM}
\end{eqnarray}
\begin{eqnarray}
\mbox{where }&& \boldsymbol{W} = \boldsymbol\Psi_n^{-1} \frac{1}{n} \left( \begin{array}{c}
\left[f_{i_0}(t_{i_0};({\boldsymbol{\beta}}, \phi))^\tau  K_{1i_0}(t_{i_0};({\boldsymbol{\beta}}, \phi)) - \gamma_{1i_0}\right] {\boldsymbol{x}_i} \\
f_{i_0}(t_{i_0};({\boldsymbol{\beta}}, \phi))^\tau  K_{2i_0}(t_{i_0};({\boldsymbol{\beta}}, \phi)) - \gamma_{2i_0}
\end{array} \right) \nonumber\\
&&~ 
~~-\boldsymbol{Q}({\boldsymbol{\theta}}_0)^{-1}\boldsymbol\Psi_n^{(0)}({{\boldsymbol{\theta}}_0 })^T \frac{1}{n} \left( \begin{array}{c}
f_{i_0}(t_{i_0};{\boldsymbol{\theta}}_0)^\tau  \boldsymbol{u}_{1i_0}^{(0)}(t_{i_0};{\boldsymbol{\theta}}_0) 
- \boldsymbol\gamma_{1i_0}^{(0)}  
\\
f_{i_0}(t_{i_0};{\boldsymbol{\theta}}_0)^\tau  \boldsymbol{u}_{2i_0}^{(0)}(t_{i_0};{\boldsymbol{\theta}}_0) 
- \boldsymbol\gamma_{2i_0}^{(0)}
\end{array} \right),~~~~\nonumber
\end{eqnarray}
with $K_{ji_0}(t_{i_0};({\boldsymbol{\beta}}, \phi))$ as defined in Section S3 of the Online Supplement for $j=1,2$, 
$\boldsymbol{u}_{1i}^{(0)}(y_i; ({\boldsymbol{\beta}},\phi))$ and  $\boldsymbol{u}_{2i}^{(0)}(y_i; ({\boldsymbol{\beta}},\phi))$  
denoting the restricted derivative of 
$\log~f_i(y_i; ({\boldsymbol{\beta}},\phi))$ with respect to ${\boldsymbol{\beta}}$ and $\phi$ under $H_0$, and 
$\boldsymbol\Psi_n^{(0)}$ being the matrix $\boldsymbol\Psi_n$ constructed using $(\boldsymbol{u}_{1i}^{(0)}, \boldsymbol{u}_{1i}^{(0)})$ 
in place of ${\boldsymbol{u}_{i}} = (\boldsymbol{u}_{1i},~\boldsymbol{u}_{2i})^T$.


\noindent
\textbf{Example \ref{SEC:glm}.1.}
Consider testing the first $r$ components ($r\leq p$) of the regression coefficient ${\boldsymbol{\beta}}$ at 
a pre-fixed value ${\boldsymbol{\beta}}_0^{(1)}$, 
the null hypothesis given by (\ref{EQ:9Gen_lin_hypothesis}) with 
$\boldsymbol{L}=\begin{bmatrix}
    ~~\boldsymbol{I}_r~ \\
    ~~\boldsymbol{O}_{(p-r)\times r}~
\end{bmatrix}$. 

We partition the relevant vectors and matrices as 
${\boldsymbol{\beta}}=({\boldsymbol{\beta}}_0^{(1)}, ~ {\boldsymbol{\beta}}_0^{(2)})$, ${\boldsymbol{x}_i}=({\boldsymbol{x}_i}^{(1)}, ~ {\boldsymbol{x}_i}^{(2)})$ and $\boldsymbol{X}=[\boldsymbol{X}_1 ~ \boldsymbol{X}_2]$, 
where ${\boldsymbol{\beta}}_0^{(1)}$ and ${\boldsymbol{x}_i}^{(1)}$ are $r$-vectors and $\boldsymbol{X}_1$ is the $n\times r$ matrix 
consisting of the first $r$ columns of $\boldsymbol{X}$. 
Let
\begin{equation}
\boldsymbol{J}_{11}=\begin{bmatrix}
     ~~\boldsymbol{J}_{11}^{11} & \boldsymbol{J}_{11}^{12}~ \\
     ~~(\boldsymbol{J}_{11}^{12})^T & \boldsymbol{J}_{11}^{22}~
\end{bmatrix},
~
\boldsymbol{V}_{11}=\begin{bmatrix}
     ~~\boldsymbol{V}_{11}^{11} & \boldsymbol{V}_{11}^{12}~ \\
     ~~(\boldsymbol{V}_{11}^{12})^T & \boldsymbol{V}_{11}^{22}~
\end{bmatrix},
~
\boldsymbol{J}_{11.2}^{-1}=\begin{bmatrix}
     ~~\boldsymbol{J}_{11.2}^{-11} & \boldsymbol{J}_{11.2}^{-12}~ \\
     ~~(\boldsymbol{J}_{11.2}^{-12})^T & \boldsymbol{J}_{11.2}^{-22}~
\end{bmatrix},\nonumber
\end{equation}
where the first block of each partitioned matrix is of order $r\times r$.

Here the asymptotic distribution of the DPD-based test statistics 
$S_{\gamma}( {{\boldsymbol{\theta}}_n^\tau}, {\widetilde{{\boldsymbol{\theta}}}_n^\tau}) $ under the null is 
that of $\sum_{i=1}^r ~  \zeta_i^{\gamma, \tau}({\boldsymbol{\theta}}_0)Z_i^2,$
where $Z_1, \cdots,Z_r$ are independent standard normal variables, and
$\zeta_1^{\gamma, \tau}({\boldsymbol{\theta}}_0), \cdots, \zeta_r^{\gamma, \tau}({\boldsymbol{\theta}}_0)$ are 
$r$ nonzero eigenvalues of the matrix 
$(1+\gamma)\boldsymbol{J}_{11,\gamma}^{11}\boldsymbol{J}_{11.2}^{-11}\boldsymbol{V}_{11}^{11}\boldsymbol{J}_{11.2}^{-11}$.
We have 
\begin{eqnarray}
\boldsymbol{W} &=& \boldsymbol\Psi_n^{-1} \frac{1}{n} \left( \begin{array}{c}
\boldsymbol{0}_r \\
\left[f_{i_0}(t_{i_0};({\boldsymbol{\beta}}, \phi))^\tau  K_{1i_0}(t_{i_0};({\boldsymbol{\beta}}, \phi)) 
- \gamma_{1i_0}\right]{\boldsymbol{x}_i}^{(2)} 
\\
f_{i_0}(t_{i_0};({\boldsymbol{\beta}}, \phi))^\tau  K_{2i_0}(t_{i_0};({\boldsymbol{\beta}}, \phi)) - \gamma_{2i_0}
\end{array} \right). \nonumber
\end{eqnarray}
The second order IF follows from (\ref{EQ:10IF_DPD_test_GLM}) with this form of $\boldsymbol{W}$.
As expected, there is no influence of contamination on the first $r$ components of the RMDPDE.
\qed


\section{Numerical Illustration: A Data Example}\label{SEC:Real_data}



We consider the multiple regression model using the 
``Salinity data" (\citet[][Chapter 2]{Rousseeuw/Leroy:1987}).
These data were discussed in \cite{Ruppert/Carroll:1980}, \cite{Rousseeuw/Leroy:1987}, and \cite{Ghosh/Basu:2013}. 
Our analysis shows, except for two potential outliers,  
that the data is well modeled by a multiple linear regression model,
taking salinity as the response variable and 
the covariates as  salinity in two weeks lag ($x_1$), 
the number of biweekly periods elapsed since the beginning of spring ($x_2$),
and the volume of river discharge into the sound ($x_3$). 
Cases 5 and 16 are outlying observations 
that correspond to periods of very heavy discharge.

The maximum likelihood estimate of the regression coefficient ${\boldsymbol{\beta}}=(\beta_0, \beta_1, \beta_2, \beta_3)^T$ 
and the error standard deviation $\sigma$ for the full data are $(9.6, 0.8, -0.03, -0.3)^T$ and $1.23$.
After deleting the outlying observations these estimates are
$(23.39, 0.70, -0.25, -0.84)^T$ and $0.91$, respectively, indicating the dramatic effect of outliers. 
\cite{Ghosh/Basu:2013} showed that the MDPDE with $\tau\geq 0.25$ can successfully generate robust estimators 
even under presence of the two outlying observations. 
In particular, the MDPDEs at $\tau=0.5$ and $\tau=1$ are, respectively, 
$\widehat{\boldsymbol{\beta}}=(18.4, 0.72, -0.2, -0.63)^T$, $\widehat\sigma=0.87$, 
and $\widehat{\boldsymbol{\beta}}=(19.19, 0.71, -0.18, -0.66)^T$, $\widehat\sigma=0.87$. 
These estimates are quite close to the outlier deleted MLE.

\begin{figure}
	\centering
	\subfloat[$H_0 : {\boldsymbol{\beta}}$ = (19.19, 0.71, $-0.18$, $-0.66$) ($\sigma = 1.23$ known)]{
		\includegraphics[width=0.32\textwidth]{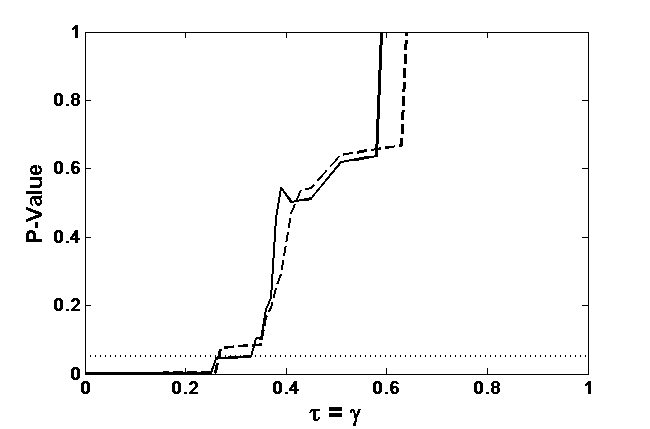}
		\label{FIG:9Salinity_k11}}
	~ 
	\subfloat[$H_0 : {\boldsymbol{\beta}}$ = (19.19, 0.71, $-0.18$, $-0.66$) ($\sigma = 0.71$ known)]{
		\includegraphics[width=0.32\textwidth]{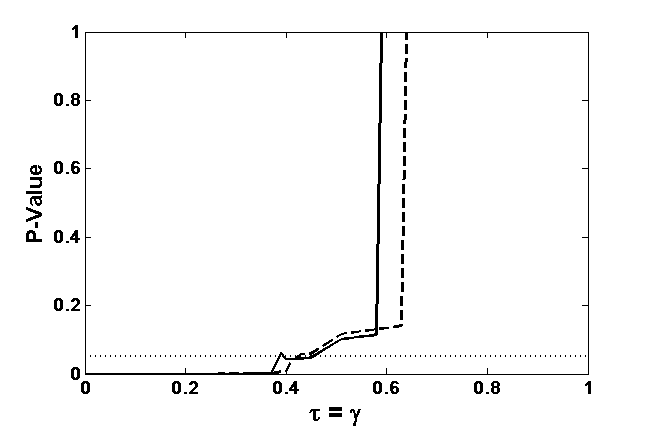}
		\label{FIG:9Salinity_k12}}
	~ 
	\subfloat[$H_0 : {\boldsymbol{\beta}}$ = (19.19, 0.71, $-0.18$, $-0.66$) ($\sigma$ unknown)]{
		\includegraphics[width=0.32\textwidth]{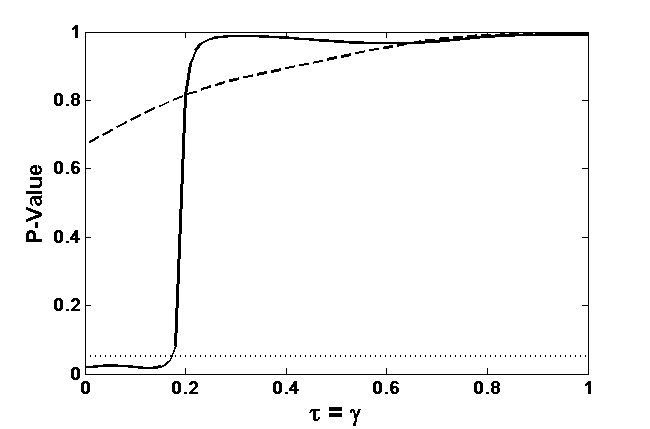}
		\label{FIG:9Salinity_1}}
	\\
	\subfloat[$H_0 : {\boldsymbol{\beta}}$ = (18.4, 0.72, $-0.2$, $-0.63$) ($\sigma = 1.23$ known)]{
		\includegraphics[width=0.32\textwidth]{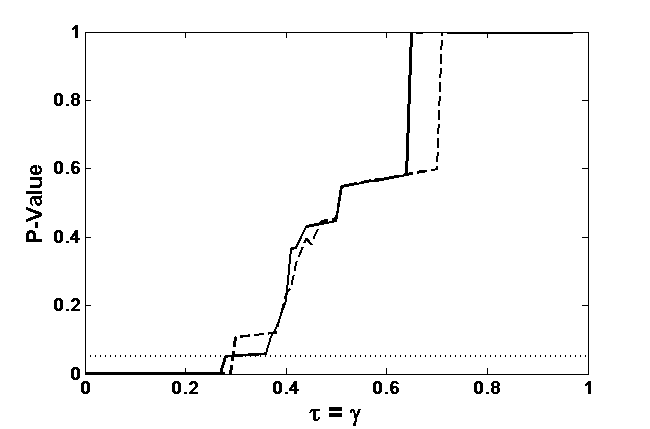}
		\label{FIG:9Salinity_k21}}
	~ 
	\subfloat[$H_0 : {\boldsymbol{\beta}}$ = (18.4, 0.72, $-0.2$, $-0.63$) ($\sigma = 0.71$ known)]{
		\includegraphics[width=0.32\textwidth]{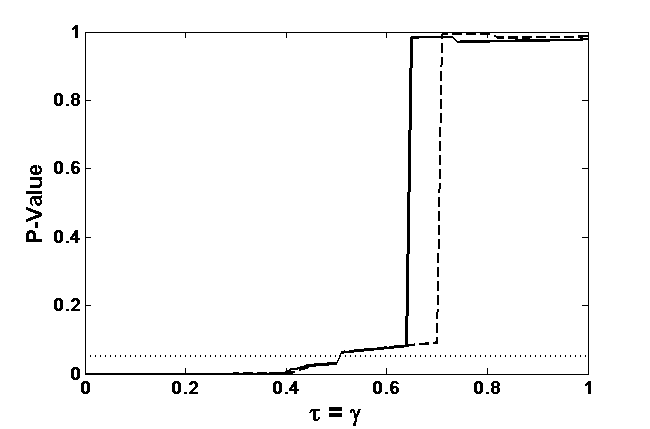}
		\label{FIG:9Salinity_k22}}
	~ 
	\subfloat[$H_0 : {\boldsymbol{\beta}}$ = (18.4, 0.72, $-0.2$, $-0.63$) ($\sigma$ unknown)]{
		\includegraphics[width=0.32\textwidth]{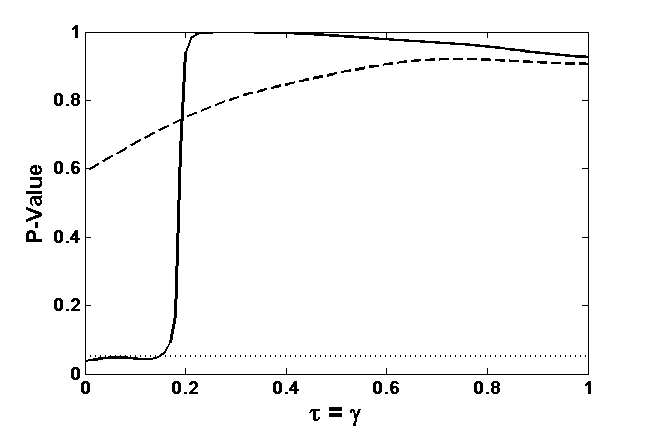}
		\label{FIG:9Salinity_2}}
	\\
	\subfloat[$H_0 : {\boldsymbol{\beta}}$ = (9.6, 0.8, $-0.03$, $-0.3$) ($\sigma = 1.23$ known)]{
		\includegraphics[width=0.32\textwidth]{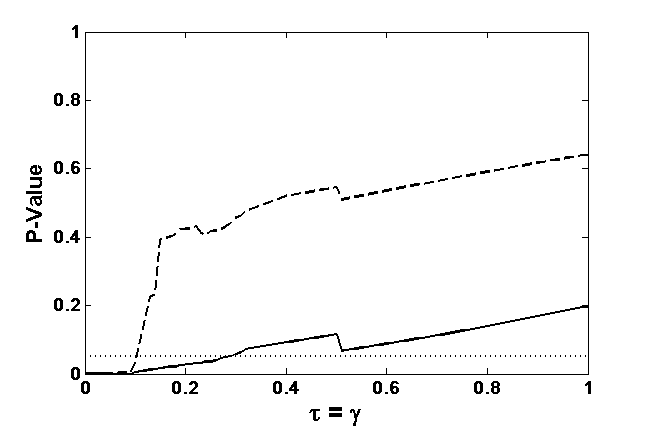}
		\label{FIG:9Salinity_k31}}
	~ 
	\subfloat[$H_0 : {\boldsymbol{\beta}}$ = (9.6, 0.8, $-0.03$, $-0.3$) ($\sigma = 0.71$ known)]{
		\includegraphics[width=0.32\textwidth]{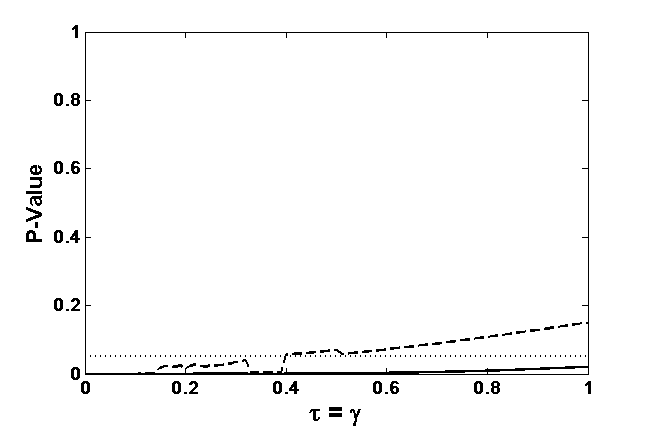}
		\label{FIG:9Salinity_k32}}
	~ 
	\subfloat[$H_0 : {\boldsymbol{\beta}}$ = (9.6, 0.8, $-0.03$, $-0.3$) ($\sigma$ unknown)]{
		\includegraphics[width=0.32\textwidth]{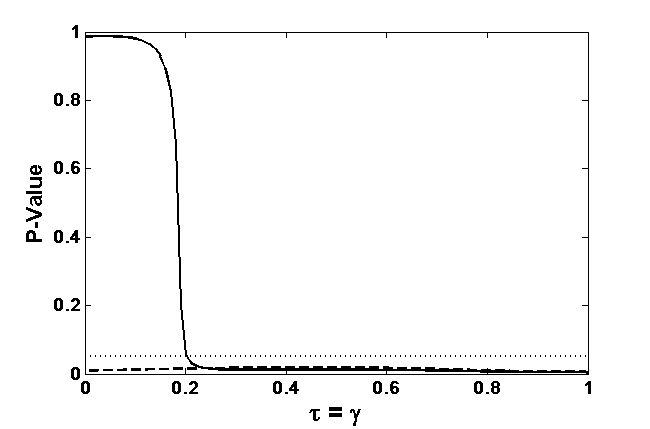}
		\label{FIG:9Salinity_3}}
	\caption{The p-values of the DPD-based tests for different $H_0$, 
	with known and unknown $\sigma^2$, for the Salinity data (Solid line - full data; dashed line - outlier deleted data)}
	\label{FIG:9Salinity_pValue}
\end{figure}

We applied the proposed DPD-based test using the full and outlier deleted data. 
We tested such hypotheses on ${\boldsymbol{\beta}}$ as 
$H_0 : {\boldsymbol{\beta}} = (19.19,~ 0.71,~ -0.18,~ -0.66)^T$, 
$H_0 : {\boldsymbol{\beta}} = (18.4, 0.72, -0.2, -0.63)^T$,
and $H_0 : {\boldsymbol{\beta}} = (9.6, 0.8, -0.03, -0.3)^T$.
They were chosen at the estimated values for two robust estimators, 
MDPDE at $\tau =1$ and $0.5$, and the non-robust MLE, respectively. 
Therefore, a robust test should accept the first two hypotheses while rejecting the third. 
We considered both simple and composite tests
by assuming $\sigma$ to be known and unknown. 
For the known $\sigma$ case, we assumed two distinct values of $\sigma$: 
$1.23$ (a non-robust estimate, MLE) and $0.71$ (a robust estimate, MDPDE at $\tau=1$).
The p-values of the proposed DPD based tests for these cases 
are presented in Figure  \ref{FIG:9Salinity_pValue}.

When $\sigma$ is unknown, the DPD-based tests
with $\tau=\gamma \geq 0.2$  give robust results by failing to reject the first two hypotheses 
(Figures \ref{FIG:9Salinity_1}, \ref{FIG:9Salinity_2})
and by rejecting the third one (Figure \ref{FIG:9Salinity_3}) under full data.
The performances of the LRT at $\tau=\gamma=0$ is clearly non-robust under full data.
When $\sigma$ is robustly specified, 
under full data the DPD-based tests still fail to reject the first two hypotheses at larger $\tau=\gamma\geq 0.5$
but the LRT rejects them (Figures \ref{FIG:9Salinity_k12}, \ref{FIG:9Salinity_k22}).
All DPD-based tests, including the LRT, successfully reject the third hypothesis under full data
for correctly specified robust $\sigma$ (Figure \ref{FIG:9Salinity_k32}). 
When $\sigma$ is incorrectly specified, the DPD-based tests at $\tau=\gamma\geq 0.5$ 
still lead to robust inference while the LRT provides incorrect inference for the first two hypotheses
(Figures \ref{FIG:9Salinity_k11}, \ref{FIG:9Salinity_k21});
the third hypothesis gets accepted by the DPD-based tests at larger $\tau=\gamma$ 
due to the incorrect specification of $\sigma$ (Figure \ref{FIG:9Salinity_k31}).

\section{Conclusions}\label{SEC:conclusion}

In this paper we have presented a general framework based on density power divergence
for performing robust tests of hypothesis in the independent but non-homogeneous 
case. We have established the wide scope of the test, and numerically demonstrated its 
applicability to the linear regression problem. Due to the generality of
the method and theoretical indicators it is expected that it will be a powerful tool 
for the practitioner, although further numerical studies would be helpful to explore 
the performance of these tests in specific situations.

Among possible extensions, we hope to study the multisample problem.
%
Another extension would be to the case of heteroscedastic models.
Finally, the choice of tuning parameters requires a thorough study.

\section*{Supplementary Material}
The Online Supplement to this paper contains the details of 
S1. required assumptions;
S2.	MDPDE under the linear regression model with fixed design along with 
some additional results for testing; 
S3.	MDPDE under the GLM with fixed design;
S4. proofs;
S5.	simulations under linear regression models;
S6.	two more data applications from linear and Poisson regression models; and
S7.	comments on the choice of tuning parameters.

\noindent
\textbf{Acknowledgment:}
The authors thank the referees, the Editor and an associate editor 
for their constructive comments and suggestions.





\end{document}